\documentclass[12pt]{amsart}
\usepackage{amsmath}
\usepackage{a4wide}
\usepackage[utf8]{inputenc}
\usepackage{amssymb}
\usepackage{amsopn}
\usepackage{epsfig}
\usepackage{amsfonts}
\usepackage{latexsym}
\usepackage{amsthm}
\usepackage{enumerate}
\usepackage[UKenglish]{babel}
\usepackage{verbatim}
\usepackage{color}
\usepackage{calrsfs}
\usepackage{pgf}
\usepackage{tikz}

\usepackage{mathrsfs}   

 \usetikzlibrary{arrows,automata}

\DeclareMathAlphabet{\pazocal}{OMS}{zplm}{m}{n}
\newtheorem{theorem}{Theorem}[section]
\newtheorem{lemma}[theorem]{Lemma}
\newtheorem{proposition}[theorem]{Proposition}
\newtheorem{corollary}[theorem]{Corollary}

\theoremstyle{definition}

\newtheorem{example}[theorem]{Example}

\theoremstyle{remark}
\newtheorem{remark}[theorem]{Remark}

\numberwithin{equation}{section}

\newcommand{\R}{\ensuremath{\mathbb{R}}}
\newcommand{\N}{\ensuremath{\mathbb{N}}}

\renewcommand{\u}{\ensuremath{\mathbf u}}

\renewcommand{\v}{\ensuremath{\mathbf{v}}}

\newcommand{\w} { {\mathbf{w}}}
\newcommand{\set}[1]{\left\{#1\right\}}
\newcommand{\la}{\lambda}
\newcommand{\ga}{\gamma}
\newcommand{\Ga}{\Gamma}
\newcommand{\ep}{\varepsilon}
\newcommand{\f}{\infty}
\newcommand{\de}{\delta}

\newcommand{\lle}{\preccurlyeq}
\newcommand{\lge}{\succcurlyeq}

\newcommand{\si}{\sigma}
\newcommand{\La}{\Lambda}

\begin{document}

\title{How likely can a point be in different Cantor sets}

 \author[K.Jiang]{Kan Jiang}
\address[K. Jiang]{Department of Mathematics, Ningbo University, Ningbo, Zhejiang, People's Republic of China}
\email{kanjiangbunnik@yahoo.com}

\author[D.Kong]{Derong Kong}
\address[D. Kong]{College of Mathematics and Statistics, Chongqing University, 401331, Chongqing, P.R.China}
\email{derongkong@126.com}

\author[W.Li]{Wenxia Li}
\address[W. Li]{School of Mathematical Sciences, Shanghai Key Laboratory of PMMP, East China Normal University, Shanghai 200062,
People's Republic of China}
\email{wxli@math.ecnu.edu.cn}

\dedicatory{}



\begin{abstract}
Let $m\in\N_{\ge 2}$, and let $\mathcal K=\set{K_\lambda: \lambda\in(0, 1/m]}$ be a class of Cantor sets, where $K_{\lambda}=\set{\sum_{i=1}^\infty d_i\lambda^i: d_i\in\{0,1,\ldots, m-1\}, i\ge 1}$. We investigate in this paper the likelyhood of a fixed point in the Cantor sets of $\mathcal K$. More precisely, for a fixed point $x\in(0,1)$ we consider the parameter set $\Lambda(x)=\set{\lambda\in(0,1/m]: x\in K_\lambda}$, and show that $\Lambda(x)$ is a topological Cantor set having zero Lebesgue measure and full Hausdorff dimension. Furthermore, by constructing a sequence of Cantor subsets with large thickness in $\Lambda(x)$  we prove that the intersection $\Lambda(x)\cap\Lambda(y)$ also has full Hausdorff dimension for any $x, y\in(0,1)$.
\end{abstract}
\maketitle

\section{Introduction}\label{s1}

Given an integer $m\ge 2$, for $\la\in{(0,1/m]}$ let $K_\la$ be the {self-similar} set generated by the \emph{iterated function system} (IFS)
$
\Psi_\la:=\set{f_i(x)=\la (x+i): i=0,1,\ldots,m-1}.
$
In other words,  $K_\la$ is the unique non-empty compact set in  $\R$ satisfying  $K_\la=\bigcup_{i=0}^{m-1}f_i(K_\la)$. So we can rewrite $K_\la$ as
\begin{equation}\label{eq:self-similar-set}
K_\la=\set{\sum_{i=1}^\f d_i\la^i: d_i\in\set{0,1,\ldots, m-1}, i\ge 1}.
\end{equation}  When $\la=1/m$, it is clear that $K_{1/m}=[0,1]$. When $\la\in(0,1/m)$, it is easy to see that $K_\la$ is a Cantor set with its
  convex hull   $conv(K_\la)=[0, \frac{(m-1)\la}{1-\la}]$; and the Hausdorff dimension of $K_\la$ is given by $-\log m/\log\la$ (cf.~\cite{Hutchinson_1981}).

  Let $\mathcal K=\set{K_\la: \la\in(0,1/m]}$ be the class of all self-similar sets defined as in (\ref{eq:self-similar-set}). We are interested in the following question: \emph{how likely can a fixed point be in different self-similar sets of $\mathcal K$?}
In other words, for a point $x\in[0,1]$  the above question is to ask how large of  the parameter set
\[
\La(x):=\set{\la\in(0,1/m]: x\in K_\la}.
\]
Observe that the point  $0$ always belongs to $K_\la$ for all $\la\in(0,1/m]$. Then $\La(0)=(0,1/m]$. Furthermore, $1\in K_\la$ if and only if $\la=1/m$. This implies  $\La(1)=\set{1/m}$. So, it is interesting to study $\La(x)$   for $x\in (0,1)$. The main difficulty in the study of $\La(x)$ is that $\La(x)$ always involves infinitely many (a continuum of) self-similar sets in $\mathcal K$. We will prove in Theorem \ref{th:category-La(x)}   that $\La(x)$ is a \emph{topological Cantor set}: a nonempty compact set having neither interior nor isolated points.  And $\La(x)$ is a Lebesgue null set but has full Hausdorff dimension (see Corollary  \ref{cor:measure-dimension-La(x)}). Furthermore, by carefully constructing a sequence of Cantor subsets with large thickness in $\La(x)$ we show in Theorem \ref{th:intersection-La(x)-La(y)} that for any two points $x, y\in(0,1)$   the intersection $\La(x)\cap\La(y)$ also has full Hausdorff dimension. As a by-product we prove that $\La(x)$ contains arbitrarily long arithmetic progressions (see Remark \ref{rem:1}).

Our motivation to study  $\La(x)$ comes from unique beta expansions. When $\la\in(1/m,1)$ the IFS $\Psi_\la=\set{f_i(x)}_{i=0}^{m-1}$ has overlaps, that is, for the attractor $K_\la$ we have $f_i(K_\la)\cap f_{i+1}(K_\la)\ne\emptyset$ for all $0\le i< m-1$. Denote by $U_\la$ the set of $x\in K_\la$ such that $x$ has a unique coding, i.e., there exists a unique sequence $(d_i)\in\set{0,1,\ldots, m-1}^\N$ such that $x=\sum_{i=1}^\f d_i\la^i$. The last two authors and their coauthors recently studied in \cite{KLLWX-2020} the topological and fractal properties of the set $U(x):=\set{\la\in(1/m,1): x\in U_\la}$. They showed that $U(x)$ has  zero Lebesgue measure for all $x>0$, and $\dim_H U(x)=1$ for any $x\in(0,1]$. Moreover, for typical $x>0$ the set $U(x)$ contains isolated points. Dajani et al.~\cite{Dajani-Komornik-Kong-Li-2018} showed that for any $x\in(0,1]$ the algebraic sum $U(x)+U(x)$ contains interior points. When  $x=1$, the set $U(1)$ was extremely studied. Komornik and Loreti \cite{Komornik_Loreti_2007} showed that the topological closure of $U(1)$ is a  {topological Cantor set}. Allaart and the second author studied in \cite{Allaart-Kong-2018}   the local dimension of $U(1)$, and showed that $U(1)$ has more weight close to $1/m$. Bonanno et al.~\cite{Bon-Car-Ste-Giu-2013} found   its close connection to the bifurcation set of $\alpha$-continued fractions and many other dynamical systems.

Another motivation to study $\La(x)$ is from the work of Boes, Darst and Erd\H os   \cite{BDE-1981}, where they considered a similar class of Cantor sets $C_\la$ in $[0,1]$ (in fact, they considered a class of fat Moran sets with positive Lebesgue measures). They showed that
$\La'(x):=\set{\la\in(0,1/2): x\in C_\la}$
 is of first category for any $x\in(0,1)$, and then they concluded that for a  second category set of $\la\in(0,1/2)$ the Cantor set $C_\la$ excluding the two endpoints   contains only irrational numbers.
\begin{figure}[h!]
\begin{center}
\begin{tikzpicture}[
    scale=40,
    axis/.style={very thick, ->},
    important line/.style={thick},
    dashed line/.style={dashed, thin},
    pile/.style={thick, ->, >=stealth', shorten <=2pt, shorten
    >=2pt},
    every node/.style={color=black}
    ]

    \draw[important line] ({1/3}, 0)--({1/2}, 0);
     \node[] at({1/3}, 0.01){$\frac{1}{3}$};  \node[] at({1/2}, 0.01){$\frac{1}{2}$};

     \draw[important line] ({1/3}, -0.03)--({0.366025}, -0.03);  \draw[important line] ({0.396608}, -0.03)--({1/2}, -0.03);
     \node[] at({0.366025}, -0.02){$110^\f$};  \node[] at({0.396608}, -0.02){$101^\f$};

 \draw[important line] ({1/3}, -0.06)--(0.342508, -0.06);\draw[important line](0.352201,-0.06)--({0.366025}, -0.06);
  \draw[important line] ({0.396608}, -0.06)--(0.423854, -0.06);\draw[important line] ({0.435958}, -0.06)--({1/2}, -0.06);

 \draw[important line] ({1/3}, -0.09)--(0.336197, -0.09); \draw[important line] (0.339163, -0.09)--(0.342508, -0.09);
 \draw[important line](0.352201,-0.09)--(0.356635, -0.09);\draw[important line](0.36051,-0.09)--({0.366025}, -0.09);
  \draw[important line] ({0.396608}, -0.09)--(0.405946, -0.09); \draw[important line] (0.410811, -0.09)--(0.423854, -0.09);
  \draw[important line] ({0.435958}, -0.09)--(0.456553, -0.09);\draw[important line] (0.461249, -0.09)--({1/2}, -0.09);
\end{tikzpicture}
\end{center}
\caption{A geometrical construction of $\La(x)$ with $x=1/2$ and $m=2$. See Example \ref{ex:1} for more explanation.}\label{fig:La(x)}
\end{figure}
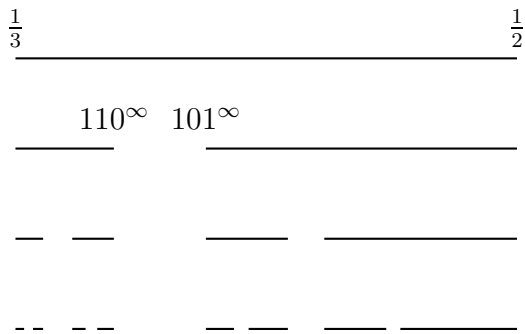

Now we state our main results. First we show that $\La(x)$ is a topologically small set.
\begin{theorem}\label{th:category-La(x)}
  For any  $x\in(0,1)$  the set $\La(x)$ is a topological Cantor set with $\min\La(x)=\frac{x}{m-1+x}$ and $\max\La(x)=1/m$.
\end{theorem}

Theorem \ref{th:category-La(x)} suggests that  $\La(x)$ can be obtained  by successively removing a sequence of open intervals from the closed interval $[\frac{x}{m-1+x}, \frac{1}{m}]$. At the end of Section \ref{sec:cateogry-La(x)} we will introduce a   geometrical  construction of $\La(x)$ (see  Figure \ref{fig:La(x)} for an example).

Observe that  the self-similar set $K_\la$  is \emph{dimensional homogeneous}, i.e., for any $y\in K_\la$  we have $\lim_{\de\to 0}\dim_H(K_\la\cap(y-\de, y+\de))=\dim_H K_\la$. This means that  $K_\la$ is distributed uniformly in the dimensional sense. Our next result states that the local dimension of $\La(x)$ {behaves differently}.
\begin{theorem}\label{th:local-dimension-La(x)}
  Let $x\in(0,1)$. Then for any $\la\in\La(x)$ we have
  \[
  \lim_{\de\to 0}\dim_H(\La(x)\cap(\la-\de, \la+\de))=\dim_H K_\la=\frac{\log m}{-\log\la}.
  \]
\end{theorem}
As a direct consequence of Theorem \ref{th:local-dimension-La(x)} we obtain that $\La(x)$ has zero Lebesgue measure and full Hausdorff dimension.
\begin{corollary}\label{cor:measure-dimension-La(x)}
  For any  $x\in(0,1)$  the set $\La(x)$ is a Lebesgue null set, and $\dim_H \La(x)=1$.
\end{corollary}
\begin{remark}
  By Corollary  \ref{cor:measure-dimension-La(x)}  it follows that  $\bigcup_{x\in(0,1)\cap\mathbb Q}\La(x)$  has zero Lebesgue measure. This  implies  that
  \[(0,1/m]\setminus\bigcup_{x\in(0,1)\cap\mathbb Q}\La(x)=\set{\la\in(0,1/m]: K_\la\subset\big((0,1)\cap\mathbb Q\big)^c}\]
  has full Lebesgue measure $1/m$. In other words, for Lebesgue almost every $\la\in(0,1/m]$ the set $K_\la\setminus\set{0}$ contains only irrational numbers.
\end{remark}
Another consequence of   Theorem \ref{th:local-dimension-La(x)} gives the following result for the segment of $\La(x)$.
\begin{corollary}\label{cor:segment-La(x)}
  \label{cor:segment-La(x)}
  Let $x\in(0,1)$. Then for any open interval    $(a,b)$ with $\La(x)\cap(a,b)\ne\emptyset$ we have
    \[
  \dim_H(\La(x)\cap (a,b))=\sup_{\la\in\La(x)\cap(a,b)}\dim_H K_\la.
  \]
\end{corollary}
Note by Theorem \ref{th:category-La(x)} that $\La(x)$ is a topological Cantor set for any $x\in(0,1)$. So by Corollary \ref{cor:segment-La(x)} it follows that the function  $h_x: \la\to \dim_H (\La(x)\cap(0,\la))$ is a non-decreasing Devil's staircase for any $x\in(0,1)$.


Our last result  shows that for any two points $x, y\in(0,1)$  there are infinitely many $\la\in(0,1/m]$ such that    $K_\la$ contains both points $x$ and $y$.  In other words, we show that $\La(x)\cap\La(y)$ contains infinitely many points.
\begin{theorem}\label{th:intersection-La(x)-La(y)}
For any $x, y\in(0,1)$ we have $\dim_H(\La(x)\cap\La(y))=1$.
\end{theorem}

The paper is organized in the following way.
In Section \ref{sec:cateogry-La(x)} we show that $\La(x)$ is a {topological} Cantor set, and prove Theorem \ref{th:category-La(x)}. In Section \ref{sec:local-dimension} we study the local dimension of $\La(x)$, and prove Theorem \ref{th:local-dimension-La(x)}, Corollaries \ref{cor:measure-dimension-La(x)} and \ref{cor:segment-La(x)}.  In Section \ref{sec:intersections-La(x)} we construct a sequence of Cantor subsets $(E_k(x))$ of $\La(x)$ such that the  thickness of $E_k(x)$ goes to infinity as $k\to\f$. Based on this we show that for any $x, y\in(0,1)$ the intersection $\La(x)\cap\La(y)$ has full Hausdorff dimension, and prove Theorem \ref{th:intersection-La(x)-La(y)}.

\section{Topological structure of $\La(x)$}\label{sec:cateogry-La(x)}
 Note that for $\la\in(0,1/m)$ the IFS $\Psi_\la$ satisfies the strong separation condition. So there exists a natural bijective map from the symbolic space $\set{0,1,\ldots, m-1}^\N$ to the {self-similar} set $K_\la$ defined by
\[
\pi_\la:\set{0,1,\ldots, m-1}^\N\to K_\la;\quad (d_i)\mapsto \sum_{i=1}^\f d_i\la^{i}.
\]
The infinite sequence $(d_i)$ is called a (unique) coding of $\pi_\la((d_i))$ with respect to the IFS $\Psi_\la$.
  Now we fix $x\in (0,1)$, and define a map from $\La(x)\setminus\set{1/m}$ to $\set{0,1,\ldots,m-1}^\N$ by
\[
\Phi_x: \La(x)\setminus\set{1/m}\to \set{0,1,\ldots,m-1}^\N;\quad \la\mapsto \pi_\la^{-1}(x).
\]
So $\Phi_x(\la)$ is the unique coding of $x\in K_\la$ with respect to the IFS $\Psi_\la$. For convenience, when $\la=1/m$ we  define
$\Phi_x(1/m)$  as the greedy $m$-adic expansion of $x$, i.e., $\Phi_x(1/m)$ is the lexicographically largest sequence $(d_i)\in\set{0,1,\ldots,m-1}^\N$ such that $x=\sum_{i=1}^\f d_i/m^i$.

Now we recall some terminology from symbolic dynamics (cf.~\cite{Lind_Marcus_1995}). In this paper we fix our \emph{alphabet} $\set{0,1,\ldots, m-1}$, and let $\set{0,1,\ldots, m-1}^\N$ denote the set of all infinite sequences $(d_i)$ with each element $d_i\in\set{0,1,\ldots,m-1}$. For a \emph{word} $\w$ we mean a finite string of digits, i.e., $\w=w_1w_2\ldots w_n$ for some $n\ge 0$. Here $n$ is called the \emph{length} of $\w$. In particular, when $n=0$ we call $\w=\epsilon$ the empty word. Let  $\set{0,1,\ldots, m-1}^*$ denote the set of all finite words. Then $\set{0,1,\ldots, m-1}^*=\bigcup_{n=0}^\f\set{0,1,\ldots,m-1}^n$, where $\set{0,1,\ldots, m-1}^n$ denotes the set of all words of length $n$. For two words $\u=u_1\ldots u_k, \v=v_1\ldots v_n$ we write $\u\v=u_1\ldots u_k v_1\ldots v_n\in\set{0,1,\ldots,m-1}^*$. Similarly, for a sequence $\mathbf a=a_1a_2\ldots \in\set{0,1,\ldots, m-1}^\N$ we write $\u\mathbf a=u_1\ldots u_k a_1 a_2\ldots$. For $j\in\N$ we denote by $\u^j=\u\cdots\u$ the $j$-fold concatenation of $\u$ with itself, and by $\u^\f=\u\u\cdots$ the period sequence with period block $\u$. For a word $\u=u_1\ldots u_k$, if $u_k<m-1$, then we write $\u^+=u_1\ldots u_{k-1}(u_k+1)$. {So} $\u^+$ is also a word. In this paper we will use \emph{lexicographical order} $\prec, \lle, \succ$ and $\lge$ between sequences in $\set{0,1,\ldots, m-1}^\N$. For example, for two sequences $(a_i), (d_i)$ we say $(a_i)\prec (b_i)$ if there exists $n\in\N$ such that $a_1\ldots a_{n-1}=b_1\ldots b_{n-1}$ and $a_n<b_n$. Also, we write $(a_i)\lle (b_i)$ if $(a_i)\prec(b_i)$ or $(a_i)=(b_i)$. 

\begin{lemma}\label{lem:lexicographical-order}
Let $x\in(0,1)$. Then  $\Phi_x$ is strictly decreasing in $\La(x)$ with respect to the lexicographical order in $\set{0,1,\ldots, m-1}^\N$.
\end{lemma}
 \begin{proof}
   Let $\la_1, \la_2\in\La(x)$ with $\la_1<\la_2$. Write $(a_i)=\Phi_x(\la_1)$ and $(b_i)=\Phi_x(\la_2)$. Suppose on the contrary that $(a_i)\lle (b_i)$. Then
   \[
   x=\sum_{i=1}^\f a_i \la_1^i\le\sum_{i=1}^\f b_i\la_1^i<\sum_{i=1}^\f b_i\la_2^i=x,
   \]
   leading to a contradiction.
 \end{proof}
As a direct consequence of Lemma \ref{lem:lexicographical-order} we can determine the extreme values of $\La(x)$.
 \begin{lemma}
   \label{lem:minimum-La(x)}
   For any $x\in(0,1)$ we have $\min\La(x)=\frac{x}{m-1+x}$ and $\max\La(x)=1/m$.
 \end{lemma}
 \begin{proof}
 Note that $K_{1/m}=[0,1]$. So by the definition of $\La(x)$ it is clear that $\max\La(x)=1/m$. On the other hand,
   by Lemma \ref{lem:lexicographical-order} it follows that the smallest element $\la=\min\La(x)$ satisfies
   $\Phi_x(\la)=(m-1)^\f$. This gives that $x=\frac{(m-1)\la}{1-\la}$, and hence, $\min\La(x)=\la=\frac{x}{m-1+x}$.
 \end{proof}
 Recall that  $K_\la$ is the self-similar set generated by the IFS $\Psi_\la=\set{f_i(x)=\la(x+i)}_{i=0}^{m-1}$. Geometrically, $K_\la$ can be constructed via a decreasing sequence of nonempty compact sets. More precisely,
\[
 K_\la=\bigcap_{n=1}^\f K_\la(n) \quad \textrm{with}\quad K_\la(n) =\bigcup_{d_1\ldots d_n\in\set{0,1,\ldots, m-1}^n}f_{d_1}\circ \cdots\circ f_{d_n}\left(\Big[0, \frac{(m-1)\la}{1-\la}\Big]\right).
\]
When $\la\in(0,1/m)$, for each $n\in\N$ the   set $K_\la(n)$ is the union of $m^n$ pairwise disjoint subintervals of equal length $\la^n\frac{(m-1)\la}{1-\la}$. Furthermore, $K_\la(n+1)\subset K_\la(n)$ for any $n\ge 1$.

Denote by $d_H$ the Hausdorff metric in the space $\mathcal C(\R)$ {consisting} of all non-empty compact sets in $\R$ (cf.~\cite{Falconer_1990}). In the following lemma we show that the set sequence $\set{K_\la(n)}_{n=1}^\f$ converges to  $K_\la$ uniformly.
\begin{lemma}
  \label{lem:uniform-converges}
   $d_H(K_\la(n), K_\la)\to 0$ uniformly as $n\to \f$ for all $\la\in(0, 1/m]$.
\end{lemma}
\begin{proof}
Note that $K_\la\subset K_\la(n)$ for any $\la\in(0,1/m]$ and $n\in\N$, and $K_\la(n)$ is the union of $m^n$ pairwise disjoint subintervals of equal length $\la^n\frac{(m-1)\la}{1-\la}$. This implies
  that for any $\la\in(0,1/m]$ we have
  \[
  d_H(K_\la(n), K_\la)\le \la^n\frac{(m-1)\la}{1-\la}\le  \frac{1}{m^n}\;\to 0 \quad\textrm{as }n\to\f.
  \]
\end{proof}
 Next we show that $\La(x)$ is closed.
 \begin{lemma}
   \label{lem:closed}
For any $x\in(0,1)$ the set $\La(x)$ is closed.
 \end{lemma}
 \begin{proof}
  Suppose $(\la_j)$ is a sequence in $\La(x)$ with $\lim_{j\to\f}\la_j=\la_0$. We will show that $\la_0\in\La(x)$, i.e., $x\in K_{\la_0}$. By Lemma \ref{lem:uniform-converges} it follows that for each $\ep>0$ there exists a large  $N\in\N$ such that for any $\la\in (0,1/m]$ and any $n>N$ we have
   \begin{equation}
     \label{eq:ko30-1}
     d_H(K_\la(n), K_\la)< \frac{\ep}{3}.
   \end{equation}
   Take $n>N$. Since $\la_j\to\la_0$ as $j\to\f$, there exists $J\in\N$ such that for any $j>J$ we have
   \begin{equation}
     \label{eq:ko30-2}
     d_H(K_{\la_j}(n), K_{\la_0}(n))<\frac{\ep}{3}.
   \end{equation}
  Then  by (\ref{eq:ko30-1}) and (\ref{eq:ko30-2}) it follows that for any $j>J$ and $n>N$,
   \begin{align*}
      d_H(K_{\la_j}, K_{\la_0})
     &\le d_H(K_{\la_j}, K_{\la_j}(n))+d_H(K_{\la_j}(n), K_{\la_0}(n))+d_H(K_{\la_0}(n), K_{\la_0}) \\
     &<\frac{\ep}{3}+\frac{\ep}{3}+\frac{\ep}{3} =\ep.
   \end{align*}
   Since $\ep>0$ was arbitrary, this implies that
   \begin{equation}
     \label{eq:15-11}
     d_H(K_{\la_j}, K_{\la_0})\to 0 \quad\textrm{as }j\to\f.
   \end{equation}

   Now suppose on the contrary that $x\notin K_{\la_0}$. Note that $K_{\la_0}$ is compact. Then $dist(x, K_{\la_0})=\inf\set{|x-y|: y\in K_{\la_0}}>0$, and thus by (\ref{eq:15-11}) it follows that $dist(x, K_{\la_{j}})>0$ for  all $j\in\N$ sufficiently large. This   leads to a contradiction with $x\in K_{\la_j}$ for all $j\ge 1$. Hence,   $x\in K_{\la_0}$, and thus $\La(x)$ is closed.
  \end{proof}
To prove Theorem \ref{th:category-La(x)} we still need the following lemma.
\begin{lemma}
  \label{lem:lip-upper}
 { Let $x\in (0,1)$ and $q\in(\frac{x}{m-1+x},1/m)$.} Then there exists a constant $C>0$ such that for any two $\la_1, \la_2\in\La(x)\cap(0,q]$ we have
  \[
 |\pi_q(\Phi_x(\la_1))-\pi_q(\Phi_x(\la_2))|>C|\la_1-\la_2|.
  \]
\end{lemma}
\begin{proof}
{Note by Lemma \ref{lem:minimum-La(x)} that $\min\La(x)=\frac{x}{m-1+x}$. So $\La(x)\cap(0,q]\ne\emptyset$ for any $q>\frac{x}{m-1+x}$.}
  Take $\la_1, \la_2\in\La(x)\cap(0,q]$ with $\la_1<\la_2$. By Lemma \ref{lem:lexicographical-order} we have
  \[
  (a_i):=\Phi_x(\la_1)\succ\Phi_x(\la_2)=:(b_i).
  \]
  So there exists $n\in\N$ such that $a_1\ldots a_{n-1}=b_1\ldots b_{n-1}$ and $a_n>b_n$. Then by using $\la_1, \la_2\in\La(x)$ we obtain that
  \[
 \sum_{i=1}^{n-1}b_i\la_2^{i}\le\sum_{i=1}^\f b_i\la_2^{i}=x= \sum_{i=1}^\f a_i\la_1^i\le  \sum_{i=1}^{n-1}a_i\la_1^i+ \sum_{i=n}^\f(m-1)\la_1^i,
  \]
  where for $n=1$ we set the sum for the index over an empty set as zero. This implies that
  \begin{equation}\label{eq:29-0}
  \begin{split}
    \frac{x}{\la_1\la_2}(\la_2-\la_1)&=\frac{x}{\la_1}-\frac{x}{\la_2}\\
    &\le\sum_{i=1}^{n-1}a_i\la_1^{i-1}-\sum_{i=1}^{n-1}b_i\la_2^{i-1}+ \sum_{i=n}^\f(m-1)\la_1^{i-1}\\
    &=\sum_{i=1}^{n-1}a_i\la_1^{i-1}-\sum_{i=1}^{n-1}a_i\la_2^{i-1}+\frac{ m-1 }{ 1-\la_1 }\la_1^{n-1}\\
    &\le \frac{ m-1 }{ 1-\la_1 }\la_1^{n-1}.
  \end{split}
  \end{equation}
  Observe that $(a_i)$ and $(b_i)$ are distinct  at the first place $n$. Then
  \begin{equation}\label{eq:29-1}
  |\pi_q((a_i))-\pi_q((b_i))|\ge {q^{n}-\frac{m-1}{1-q}q^{n+1}}=\frac{1-mq}{1-q}q^n.
  \end{equation}
  Hence, by (\ref{eq:29-0}) and (\ref{eq:29-1}) we conclude that
  \begin{align*}
    |\pi_q(\Phi_x(\la_1))-\pi_q(\Phi_x(\la_2))|&=|\pi_q((a_i))-\pi_q((b_i))|\\
    &\ge\frac{1-mq}{1-q}q^n\ge\frac{1-mq}{1-q}\la_1^n\\
    &\ge\frac{1-mq}{1-q}\times\frac{(1-\la_1)x}{(m-1)\la_2}|\la_1-\la_2|\\
    &\ge\frac{(1-mq)x}{(m-1)q}|\la_1-\la_2|
  \end{align*}
  as desired.
\end{proof}
 \begin{proof}
   [Proof of Theorem \ref{th:category-La(x)}]
   By Lemmas  \ref{lem:minimum-La(x)} and \ref{lem:closed} we only need to prove that $\La(x)$ has neither isolated nor interior points. First we prove that $\La(x)$ has no isolated points. Let $\la\in\La(x)$. We split the proof into the following   two cases.

   Case I.   $\la\in(0,1/m)$. Then $\Phi_x(\la)=(d_i)$ is the unique coding of $x$ with respect to the IFS $\Psi_\la$. For $n\ge 1$ let $p_n$ be the unique root in $(0,1)$ of the equation
   \[
   x=\sum_{i=1}^{n-1}d_i p_n^i+d_n' p_n^n,
   \]
   where $d_n'=d_n+1\pmod m$.    Then for all sufficiently large $n\in\N$ we have $p_n\in(0,1/m)$ and $\Phi_x(p_n)=d_1\ldots d_{n-1}d_n' 0^\f$. So $p_n\in\La(x)$ for all large $n\in\N$. Furthermore, by Lemma \ref{lem:lip-upper} it follows that $p_n\to\la$ as $n\to\f$. Thus, $\la$ is not isolated in $\La(x)$.

   Case II. $\la=1/m$. Then $\Phi_x(1/m)=(d_i)$ is the greedy $m$-adic expansion of $x$. Since $x\in(0,1)$, there {exist} infinitely many $n\in\N$ such that $d_n<m-1$. For all such $n$ we define $p_n\in(0,1)$ so that
   $x=\sum_{i=1}^{n-1}d_i p_n^i+(m-1)p_n^n$. Then $p_n\in(0,1/m)$ and $\Phi_x(p_n)=d_1\ldots d_{n-1}(m-1)0^\f$. So $p_n\in\La(x)$ for all {$n\in\N$ satisfying} $d_n<m-1$. Again, by Lemma \ref{lem:lip-upper} it follows that $p_n$ converges to $1/m$ along a suitable subsequence, and thus $1/m$ is not isolated in $\La(x)$.

   By Cases I and II we conclude that $\La(x)$ has no isolated points.
  Next we prove that $\La(x)$ has no interior points. {It suffices to prove} that for any two points $\la_1, \la_2\in\La(x){\setminus\set{1/m}}$  there must exist $\la$ between $\la_1$ and $\la_2$ but not in $\La(x)$.

 Take $\la_1, \la_2\in\La(x){\setminus\set{1/m}}$ with $\la_1<\la_2$. Then by Lemma \ref{lem:lexicographical-order} we have
 \[
 (a_i)=\Phi_x(\la_1)\succ\Phi_x(\la_2)=(b_i).
 \]
 So there exists $n\in\N$ such that $a_1\ldots a_{n-1}=b_1\ldots b_{n-1}$ and $a_n>b_n$. We define two sequences
 \[
 \xi=b_1\ldots b_n(m-1)^\f,\quad \zeta=b_1\ldots b_{n-1}(b_n+1)0^\f.
 \]
 Then $(b_i)\lle \xi\prec \zeta\lle (a_i)$. {So by using $0<\la_1<\la_2<1/m$ it follows} that
 \begin{equation}\label{eq:30-3}
 \pi_{\la_1}(\xi)<\pi_{\la_1}(\zeta)\le\pi_{\la_1}((a_i))=x=\pi_{\la_2}((b_i))\le\pi_{\la_2}(\xi)<\pi_{\la_2}(\zeta).
 \end{equation}
 Let $I_{\la}:=(\pi_\la(\xi), \pi_\la(\zeta))$. Then (\ref{eq:30-3}) implies that $x$ is between two disjoint open intervals $I_{\la_1}$ and $I_{\la_2}$. Observe that the map $\la\mapsto \overline{I_\la}=[\pi_\la(\xi), \pi_\la(\zeta)]$ is continuous with respect to the Hausdorff metric $d_H$. So there must exist $\la\in(\la_1,\la_2)$ such that $x\in I_\la$. Since $I_\la\cap K_\la=\emptyset$, it follows that $x\notin K_\la$, i.e., $\la\in(\la_1,\la_2)\setminus\La(x)$. This completes the proof.
 \end{proof}

{ At the end of this section we describe the
geometrical construction of $\La(x)$.}  By Theorem \ref{th:category-La(x)} it follows that $\La(x)$ is a Cantor set in $\R$ for any $x\in(0,1)$. So, it can be obtained  by successively removing a sequence of open intervals  from $conv(\La(x))=[\frac{x}{m-1+x},  \frac{1}{m}]$.  Let $(x_i)=\Phi_x(1/m)$ be the greedy $m$-adic expansion of $x$. Since $x\in(0,1)$, there exists a smallest $\ell\in\N$ such that $x_\ell<m-1$. Note that $\Phi_x(\frac{x}{m-1+x})=(m-1)^\f$. We call a word $\w\in\set{0,1,\ldots, m-1}^*$ \emph{admissible} in $\Phi_x(\La(x))$ if
 \[
 (x_i)\prec\w  0^\f\prec \w(m-1)^\f\prec (m-1)^\f.
 \]
 Since $x_i=m-1$ for all $i<\ell$, it follows that any admissible word $\w$ has  length at least $\ell$. For each admissible word $\w$ we define the associated basic interval  $J_\w=[p_\w, q_\w]$  by
 \[
 \Phi_x(p_\w)=\w(m-1)^\f\quad\textrm{and}\quad \Phi_x(q_\w)=\w  0^\f.
 \]
 Then for two admissible words $\u, \w$, if $\u$ is  a prefix of $\w$, then $J_\w\subset J_\u$. Furthermore, for any admissible word $\w$ there must exist an admissible word $\v$ such that $\w$ is a prefix of $\v$, i.e., $\w$ has an offspring $\v$.  Denote by $\mathcal A(x)=\bigcup_{n=\ell}^\f A_n(x)$ the set of all admissible words, where $A_n(x)$ consists of all admissible words of length $n$. So the basic intervals $\set{J_\w: \w\in\mathcal A(x)}$ have a tree structure, and the Cantor set $\La(x)$ can be constructed geometrically as
 \[
 \La(x)=\bigcap_{n=\ell}^\f\bigcup_{\w\in A_n(x)}J_\w.
 \]

 \begin{example}
   \label{ex:1}
   Let $m=2$ and $x=1/2$. Then  $conv(\La(x))=[\frac{1}{3}, \frac{1}{2}]$. Furthermore,
   $\Phi_x(\frac{1}{3})=1^\f$ and $\Phi_x(\frac{1}{2})=10^\f$. Then any admissible word has length at least $\ell=2$. By the definition of admissible words it follows that
   \[
   A_2(x)=\set{10, 11},\quad A_3{(x)}=\set{100, 101, 110, 111},\ldots,
   \]
   and in general, for any $n\in\N$ we have
   \[A_{n+1}{(x)}=\set{1\u: \u\in\set{0,1}^n}.\]
  So, in the first step we remove the open interval $H_1=(0.366025, 0.396608)\sim(110^\f, 101^\f)$ from the convex hull $[\frac{1}{3}, \frac{1}{2}]$; and in the next step we remove two open intervals
  \begin{align*}
  H_2&=(0.342508, 0.352201)\sim(1110^\f, 1101^\f),\\
   H_3&=(0.423854, 0.435958)\sim(1010^\f, 1001^\f).
  \end{align*}
   This procedure can be continued, and after finitely many steps we can get a good approximation  of $\La(x)$ (see Figure \ref{fig:La(x)}).
 \end{example}

 \section{Fractal properties of $\La(x)$: local dimension}\label{sec:local-dimension}
In this section we will investigate the local dimension of $\La(x)$, and prove Theorem \ref{th:local-dimension-La(x)}. Our proof will be split into the following two cases: (I) local dimension of $\La(x)$ at $\la=1/m$; (II) local dimension of $\La(x)$ at $\la\in(0,1/m)$.

\subsection{Local dimension of $\La(x)$ at $\la=1/m$} Observe that $1/m\in\La(x)$ for any $x\in(0,1)$. In this part we will show that the local dimension of $\La(x)$ at $1/m$ is one.
\begin{proposition}
\label{prop:local-dim-La(x)=1/m}
For any $x\in(0,1)$ we have
\begin{equation*}\label{eq:06-1}
\lim_{\de\to 0}\dim_H\left(\La(x)\cap\Big(\frac{1}{m}-\de, \frac{1}{m}+\de\Big)\right)=1=\dim_H K_{1/m}.
\end{equation*}
\end{proposition}
Our strategy to prove Proposition \ref{prop:local-dim-La(x)=1/m} is to construct a large subset of $\La(x)\cap(1/m-\de,1/m+\de)$ with its Hausdorff dimension arbitrarily close to one.
Let $x\in(0,1)$. For $k\in\N$ we set
\[
\La_k(x):=\set{\la\in \La(x): \exists N \textrm{ such that } \si^N(\Phi_x(\la))\textrm{ does not contain $k$ consecutive zeros}},
\]
where $\si$ is the left-shift map on $\set{0,1,\ldots,m-1}^\N$. Note that for any $\la\in\La(x)\setminus\bigcup_{k=1}^\f\La_k(x)$ the coding $\Phi_x(\la)$ must end with $0^\f$.
Thus     the difference between $\La(x)$ and $\bigcup_{k=1}^\f \La_k(x)$ is at most countable. So, by the countable stability of Hausdorff dimension it follows that
\[\dim_H\La(x)=\dim_H\bigcup_{k=1}^\f\La_k(x)=\sup_{k\ge 1}\dim_H\La_k(x).\]
In the following we will give a lower bound for the Hausdorff dimension of $\La_k(x)$. Let $(x_i)=\Phi_x(1/m)$ be the greedy $m$-adic expansion of $x$. Since $x\in(0,1)$, there exist  infinitely many digits $x_i<m-1$. So there exists a subsequence $(n_i)\subset\N$ such that $x_{n_i}<m-1$ for all $i\ge 1$. For a large integer $j$, let $\ga_j$ be the unique root in $(0,1/m)$ of the equation
\[
x=\sum_{i=1}^{n_j-1}x_i \ga_j^i+(x_{n_j}+1)\ga_j^{n_j}+\sum_{i=n_j+1}^\f(m-1)\ga_j^i.
\]
Then $\Phi_x(\ga_j)=x_1\ldots x_{n_j}^+(m-1)^\f$, {and} $\ga_j\nearrow 1/m$ as $j\to\f$. Define
\[
\Ga_{k,j}(x):=\set{x_1\ldots x_{n_j}^+ d_1d_2\ldots\in\set{0,1,\ldots,m-1}^\N: (d_i)\textrm{ does not contain $k$ consecutive zeros}}.
\]
\begin{lemma}
  \label{lem:large-subset}
  Let $x\in(0,1)$ and $k\in\N$. Then for a large $j\in\N$ we have
  \[
  \Ga_{k,j}(x)\subseteq\Phi_x(\La_k(x)\cap[\ga_j,1/m)).
  \]
\end{lemma}
\begin{proof}
  Take a sequence $(c_i)\in\Ga_{k,j}(x)$. Then the equation
  $x=\sum_{i=1}^\f c_i\la^i$ determines a unique $\la\in(0,1/m)$, i.e., $\Phi_x(\la)=(c_i)$. Note that
  \[
 \Phi_x(1/m)=(x_i)\prec \Phi_x(\la)=(c_i)\lle x_1\ldots x_{n_j}^+(m-1)^\f=\Phi_x(\ga_j).
  \]
  By Lemma \ref{lem:lexicographical-order} we conclude that $\la\in[\ga_j,1/m)$. Furthermore, by the definition of $\Ga_{k,j}(x)$ it follows that
  $\si^{n_j}((c_i))$ does not contain $k$ consecutive zeros. So, $\la\in\La_k(x)$. This completes the proof.
\end{proof}

To give a lower bound   of $\dim_H(\La_k(x)\cap[\ga_j, 1/m))$ we still need the following lemma.
\begin{lemma}
  \label{lem:lip-lower}
  Let $x\in(0,1)$ and $k\in\N$. Then for a large $j\in\N$ there exists $C>0$ such that for any $\la_1, \la_2\in\Phi_x^{-1}(\Ga_{k,j}(x))$ we have
  \[
  |\pi_{\ga_j}(\Phi_x(\la_1))-\pi_{\ga_j}(\Phi_x(\la_2))|\le C|\la_1-\la_2|.
  \]
\end{lemma}
\begin{proof}
  Let $\la_1, \la_2\in\Phi_x^{-1}(\Ga_{k,j}(x))$ with $\la_1<\la_2$. Then by Lemma \ref{lem:large-subset} we have $\la_1, \la_2\in{\La_k(x)\cap}[\ga_j,1/m)$. Furthermore, $(a_i)=\Phi_x(\la_1)$ and $(b_i)=\Phi_x(\la_2)$ have a common prefix of length at least $n_j$. Since $\la_1<\la_2$, by Lemma \ref{lem:lexicographical-order} we have $(a_i) \succ (b_i)$. So there exists $n>n_j$ such that $a_1\ldots a_{n-1}=b_1\ldots b_{n-1}$ and $a_n>b_n$. Note that $\si^{n_j}((a_i))$ does not contain $k$ consecutive zeros. Then
  \[
  \la_1^{n+k}+\sum_{i=1}^n a_i\la_1^i<\sum_{i=1}^\f a_i\la_1^i=x=\sum_{i=1}^\f b_i\la_2^i<\sum_{i=1}^n a_i \la_2^i.
  \]
  Rearranging the above equation it gives that
  \begin{align*}
    \la_1^{n+k}<\sum_{i=1}^n a_i(\la_2^i-\la_1^i)<\sum_{i=1}^\f(m-1)(\la_2^i-\la_1^i)=\frac{m-1}{(1-\la_1)(1-\la_2)}(\la_2-\la_1).
  \end{align*}
  Since $\la_1,\la_2\in[\ga_j,1/m)$, this implies that
  \begin{equation}\label{eq:29-2}
    \la_1^n<\frac{m-1}{\la_1^k(1-\la_1)(1-\la_2)}(\la_2-\la_1)< \frac{m^2}{\ga_j^k(m-1)}(\la_2-\la_1).
  \end{equation}
  Note that $(a_i)$ and $(b_i)$ {have} a common prefix of length $n-1$. It follows that
  \[|\pi_{\ga_j}((a_i))-\pi_{\ga_j}((b_i))|\le \ga_j^{n-1}\frac{(m-1)\ga_j}{1-\ga_j}=\frac{m-1}{1-\ga_j}\ga_j^n.\]
  Therefore, by (\ref{eq:29-2}) it follows that
  \begin{align*}
    |\pi_{\ga_j}(\Phi_x(\la_1))-\pi_{\ga_j}(\Phi_x(\la_2))|&=|\pi_{\ga_j}((a_i))-\pi_{\ga_j}((b_i))|\\
    &\le \frac{m-1}{1-\ga_j}\ga_j^n\le \frac{m-1}{1-\ga_j}\la_1^n\\
    &<\frac{m-1}{1-\ga_j}\times\frac{m^2}{\ga_j^k(m-1)}|\la_1-\la_2|=\frac{m^2}{(1-\ga_j)\ga_j^k}|\la_1-\la_2|
  \end{align*}
  as desired.
\end{proof}

\begin{lemma}
  \label{lem:dim-lower}
  Let $x\in(0,1)$ and $k\in\N$. Then for a large $j\in\N$ we have
  \[
  \dim_H (\La_k(x)\cap[\ga_j,1/m)) \ge \frac{(k-1)\log m+\log(m-1)}{-k\log \ga_j}.
  \]
\end{lemma}
\begin{proof}
This follows by Lemmas \ref{lem:large-subset} and \ref{lem:lip-lower}   that
\begin{align*}
\dim_H(\La_k(x)\cap[\ga_j,1/m))
&\ge \dim_H\Phi_x^{-1}(\Ga_{k,j}(x))\\
&\ge\dim_H\pi_{\ga_j}(\Ga_{k,j}(x))\ge \frac{(k-1)\log m+\log(m-1)}{-k\log \ga_j},
\end{align*}
where the last inequality follows from that $\si^{n_j}(\Ga_{k,j}(x))$ is a sub-shift of finite type with the forbidden block $0^k$.
\end{proof}

\begin{proof}
  [Proof of Proposition \ref{prop:local-dim-La(x)=1/m}]
  Note that $\bigcup_{k=1}^\f\La_k(x)\subseteq\La(x)$ and the difference set $\La(x)\setminus\bigcup_{k=1}^\f\La_k(x)$ is at most countable. Furthermore, by the definition of $\La_k(x)$ we know that $\La_{k_1}(x)\subset\La_{k_2}(x)$ for any $k_1<k_2$. So
  by Lemma \ref{lem:dim-lower} it follows that for any large $j\in\N$,
  \begin{align*}
 \dim_H(\La(x)\cap[\ga_j,1/m))
    &=\dim_H\bigcup_{k=1}^\f(\La_k(x)\cap[\ga_j,1/m))\\
    &=\lim_{k\to\f}\dim_H(\La_k(x)\cap[\ga_j,1/m))\\
    &\ge\lim_{k\to\f}\frac{(k-1)\log m+\log(m-1)}{-k\log \ga_j}=\frac{\log m}{-\log \ga_j}.
  \end{align*}
  Note that $\ga_j\to 1/m$ as   $j\to\f$. This implies that
  \begin{align*}
  \lim_{\de\to 0}\dim_H(\La(x)\cap(1/m-\de, 1/m+\de))&\ge\lim_{j\to\f}\dim_H(\La(x)\cap[\ga_j, 1/m))\\
  &\ge\lim_{j\to\f}\frac{\log m}{-\log \ga_j}=1,
  \end{align*}
  completing the proof.
\end{proof}

\subsection{Local dimension of $\La(x)$ at $\la\in\La(x)\setminus\set{1/m}$} In this part we will prove Theorem \ref{th:local-dimension-La(x)} for $\la\in\La(x)\setminus\set{1/m}$.
Fix $x\in(0,1)$ and $\la\in\La(x)\setminus\set{1/m}$. Let $(x_i)=\Phi_x(\la)\in\set{0,1,\ldots, m-1}^\N$ be the unique coding of $x\in K_\la$ with respect to the IFS $\Psi_\la$. Let  $n\in\N$ large enough such that $x_1\ldots x_n\ne 0^n$, and then  let $\beta_n$ and $\ga_n$ be the unique roots in $(0,1)$ of the following equations respectively:
\[
x=\sum_{i=1}^n x_i \beta_n^i+\sum_{i=n+1}^\f(m-1)\beta_n^i\quad\textrm{and}\quad x=\sum_{i=1}^n x_i \ga_n^i.
\]
In fact, by choosing $n\in\N$ sufficiently large we can even require that
\begin{equation}\label{eq:08-1}
0<\beta_n\le \la\le \ga_n<1/m.
\end{equation}
Then $\Phi_x(\beta_n)=x_1\ldots x_n(m-1)^\f$ and $\Phi_x(\ga_n)=x_1\ldots x_n 0^\f$. Furthermore,  $\lim_{n\to\f}\beta_n=\lim_{n\to\f}\ga_n=\la$.

Recall from the previous subsection that $\La_k(x)$ consists of all ${q}\in \La(x)$ such that the tail sequence of $\Phi_x({q})$ does not contain $k$ consecutive zeros. Now for a large $n\in\N$, we also define
\[
\Ga^\la_{k,n}(x):=\set{x_1\ldots x_n d_1d_2\ldots\in\set{0,1,\ldots, m-1}^\N: (d_i)\textrm{ does not contain $k$ consecutive zeros}}.
\]
\begin{lemma}
  \label{lem:local-1}
  Let $x\in(0,1), \la\in\La(x)\setminus\set{\frac{1}{m}}$ and $k\in\N$. Then for any $\de\in(0,\min\set{\la, {\frac{1}{m}}-\la})$ there exists a large $N\in\N$ such that
  \[
  \Ga^\la_{k,n}(x)\subset\Phi_x(\La_k(x)\cap(\la-\de, \la+\de))
  \]
  for any $n\ge N$.
\end{lemma}
\begin{proof}
By (\ref{eq:08-1}) there exists   $N_1\in\N$ such that for any $n>N_1$ we have $0<\beta_n\le \la\le \ga_n<1/m$.
Take $(c_i)\in\Ga^\la_{k,n}(x)$. Then the equation $x=\sum_{i=1}^\f c_i q^i$ determines a unique $q\in(0,1/m)$, i.e.,  $\Phi_x(q)=(c_i)$. Observe that
  \[
  \Phi_x(\ga_n)=x_1\ldots x_n0^\f\lle (c_i)=\Phi_x(q)\lle x_1\ldots x_n(m-1)^\f=\Phi_x(\beta_n).
  \]
  By Lemma \ref{lem:lexicographical-order} we have $\beta_n\le q\le \ga_n$. Note that $\lim_{n\to\f}\beta_n=\lim_{n\to\f}\ga_n=\la$. Then for $\de\in(0,\min\set{\la,1/m-\la})$ there exists an integer $N>N_1$ such that for any $n>N$ we have
  \[q\in [\beta_n, \ga_n]\subset(\la-\de, \la+\de)\subset(0,1/m).\]
   Clearly, by the definition of $\Ga_{k,n}^\la(x)$ the tail sequence $\si^n((c_i))$ does not contain $k$ consecutive zeros. This proves $q\in\La_k(x)$, and hence completes the proof.
\end{proof}
Now by Lemma \ref{lem:local-1} and  the same argument as in the proof of Lemma \ref{lem:dim-lower} we obtain the following lower bound for the dimension of $\La_k(x)\cap(\la-\de, \la+\de)$.
\begin{lemma}
  \label{lem:local-2}
  Let $x\in(0,1), \la\in\La(x)\setminus\set{{\frac{1}{m}}}$ and $k\in\N$. Then for any $\de\in(0,\min\set{\la, {\frac{1}{m}}-\la})$ we have
  \[
  \dim_H (\La_k(x)\cap(\la-\de, \la+\de))\ge\dim_H\pi_{\la-\de}(\Ga_{k,n}^{\la}(x))\ge \frac{(k-1)\log m+\log(m-1)}{-k\log(\la-\de)}.
  \]
\end{lemma}

 \begin{proof}
   [Proof of Theorem \ref{th:local-dimension-La(x)}]
   By Proposition \ref{prop:local-dim-La(x)=1/m} it
 suffices to prove the theorem for $\la\in\La(x)\setminus\set{1/m}$.  By Lemma \ref{lem:local-2} it follows that for any $\de\in(0,\min\set{\la, {\frac{1}{m}}-\la})$ we have
   \begin{align*}
     \dim_H(\La(x)\cap(\la-\de, \la+\de))&=\dim_H\bigcup_{k=1}^\f(\La_k(x)\cap(\la-\de,\la+\de))\\
     &=\lim_{k\to\f}\dim_H\La_k(x)\cap(\la-\de,\la+\de)\\
     &\ge\lim_{k\to\f}\frac{(k-1)\log m+\log(m-1)}{-k\log(\la-\de)}=\frac{\log m}{-\log(\la-\de)}.
   \end{align*}
   Letting $\de\to 0$ we obtain  that
   \begin{equation}\label{eq:kk-1}
   \lim_{\de\to 0}  \dim_H(\La(x)\cap(\la-\de, \la+\de))\ge \frac{\log m}{-\log \la}.
   \end{equation}

   On the other hand, by using Lemma \ref{lem:lip-upper} it follows that
   \begin{equation}\label{eq:kk-2}
   \begin{split}
   \dim_H(\La(x)\cap(\la-\de,\la+\de))&\le\dim_H\pi_{\la+\de}\circ\Phi_x(\La(x)\cap(\la-\de,\la+\de))\\
   &\le \dim_H K_{\la+\de}=\frac{\log m}{-\log(\la+\de)}.
   \end{split}
   \end{equation}
   Letting $\de\to 0$ in (\ref{eq:kk-2}), we conclude by (\ref{eq:kk-1}) that
   \[\lim_{\de\to 0}  \dim_H(\La(x)\cap(\la-\de, \la+\de))=\frac{\log m}{-\log \la}.\]
 \end{proof}

\begin{proof}
  [Proof of Corollary \ref{cor:measure-dimension-La(x)}]
 By Theorem \ref{th:local-dimension-La(x)} it follows that
 \[
 \dim_H\La(x)\ge \lim_{n\to\f}\dim_H(\La(x)\cap(1/m-1/2^n, 1/m+1/2^n))=\dim_H K_{1/m}=1.
 \]
 So $\La(x)$ has full Hausdorff dimension.
 In the following it suffices to prove that $\La(x)$ has zero Lebesgue measure for any $x\in(0,1)$.

 By Theorem \ref{th:local-dimension-La(x)} it follows that
  for any $\la\in\La(x)\setminus\set{1/m}$ there exists $\de_\la>0$ such that $\dim_H(\La(x)\cap(\la-\de_\la, \la+\de_\la))<1$, which implies that $\La(x)\cap(\la-\de_\la, \la+\de_\la)$ is a Lebesgue null set. Note by Theorem \ref{th:category-La(x)} that for any $\ga\in(0, 1/m)$ the segment $\La(x)\cap(0, \ga]$ is compact, and then it can be covered by a finite number of open intervals $\set{(\la_i-\de_{\la_i}, \la_i+\de_{\la_i})}_{i=1}^N$ for some $\la_i\in\La(x)$ with $1\le i\le N$. Since each set $\La(x)\cap (\la_i-\de_{\la_i}, \la_i+\de_{\la_i})$ is a Lebesgue null set,  so is  $\La(x)\cap(0,\ga]$ for any $\ga<1/m$. Let $(\ga_n)$ be a sequence in $(0,1/m)$ with $\ga_n\nearrow 1/m$ as $n\to\f$. Then we conclude that $\La(x)=\set{1/m}\cup\bigcup_{n=1}^\f(\La(x)\cap(0,\ga_n])$ has zero Lebesgue measure.
 \end{proof}

\begin{proof}[Proof of Corollary \ref{cor:segment-La(x)}]
Note that the function $D: \la\mapsto \dim_H K_\la=\frac{\log m}{-\log\la}$ is continuous and strictly increasing in $(0,1/m]$. Let $(a,b)$ be an open interval such that $(a,b)\cap \La(x)\ne\emptyset$. We consider the following two cases.

Case (I). There exists a $\la_*\in\La(x)\cap(a,b)$ such that $D(\la_*)=\sup_{\la\in\La(x)\cap(a,b)}D(\la)$. Then $\La(x)\cap(\la_*, b)=\emptyset$. Take $\ep>0$. By Theorem \ref{th:local-dimension-La(x)} it follows that
\begin{equation}\label{eq:k30-0}
\dim_H(\La(x)\cap(a,b))\ge \dim_H(\La(x)\cap(\la_*-\de, \la_*+\de))\ge D(\la_*)-\ep
\end{equation}
for $\de>0$ sufficiently small. On the other hand,  by Theorem \ref{th:local-dimension-La(x)} it follows that
  for any other $\la\in\La(x)\cap[a,b)=\La(x)\cap[a,\la_*]$ there must exist a $\de_\la>0$ such that
\begin{equation}\label{eq:k30-1}
\dim_H(\La(x)\cap(\la-\de_\la, \la+\de_\la))< D(\la)+\ep.
\end{equation}
Observe that the union of $(\la-\de_\la, \la+\de_\la)$ with $\la\in\La(x)\cap[a,b)$   is an open cover of the compact set $\La(x)\cap[a,b)=\La(x)\cap[a,\la_*]$. So there exists a subcover of $\La(x)\cap[a,b)$, say $\set{(\la_i-\de_{\la_i}, \la_i+\de_{\la_i}): 1\le i\le N}$, such that $\La(x)\cap[a,b)\subset\bigcup_{i=1}^N (\La(x)\cap(\la_i-\de_{\la_i}, \la_i+\de_{\la_i}))$.
Therefore, by (\ref{eq:k30-1}) we conclude that
\begin{align*}
\dim_H(\La(x)\cap(a,b))&=\dim_H(\La(x)\cap[a,b))=\max_{1\le i\le N}\dim_H (\La(x)\cap(\la_i-\de_{\la_i}, \la_i+\de_{\la_i}))\\
&\le\max_{1\le i\le N}(D(\la_i)+\ep)\le D(\la_*)+\ep.
\end{align*}
Since $\ep>0$ was taken arbitrarily, this together with (\ref{eq:k30-0}) implies that
\[
\dim_H(\La(x)\cap(a,b))=\dim_H(\La(x)\cap(a,\la_*])=D(\la_*)=\sup_{\la\in\La(x)\cap(a,b)}\dim_H K_\la
\]
as desired.

Case (II). There is no $\la_*\in\La(x)\cap(a,b)$ such that $D(\la_*)=\sup_{\la\in\La(x)\cap(a,b)}D(\la)$. Then there exists a sequence $(\la_n)\subset\La(x)\cap(a,b)$ such that $\la_n\nearrow b$ as $n\to\f$,  and $\sup_{n\ge 1}D(\la_n)= \sup_{\la\in\La(x)\cap(a,b)}D(\la)$. By Case (I) we know that for any $n\ge 1$ we have
\[
\dim_H (\La(x)\cap(a,\la_{n}])=D(\la_n).
\]
By the countable stability of Hausdorff dimension we conclude that
\[
\dim_H(\La(x)\cap(a,b))=\dim_H\bigcup_{n=1}^\f(\La(x)\cap(a,\la_n])=\sup_{n\ge 1}D(\la_n)=\sup_{\la\in\La(x)\cap(a,b)}D(\la).
\]
This completes the proof.
\end{proof}

\section{Large intersection of $\La(x)$ and $\La(y)$}\label{sec:intersections-La(x)}

In this section we will show  that the intersection  of any two sets $\La(x)$ and $\La(y)$ contains a Cantor subset of Hausdorff dimension arbitrarily close to one, and then prove Theorem \ref{th:intersection-La(x)-La(y)}.
Our strategy is to show that $\La(x)$ contains a Cantor subset with its thickness arbitrarily large.

\subsection{Construction of Cantor subsets of $\La(x)$ with large thickness} First we recall from \cite{Newhouse-1970} the definition of thickness for a Cantor set in $\R$. Suppose $E\subset \R$ is a Cantor set. Then $E$ can be obtained by successively removing countably many pairwise disjoint open intervals $\set{U_i}_{i=1}^\f$ from the closed interval $I=conv(E)$, where $conv(E)$ denotes the convex hull of $E$. In the first step, we remove the open interval $U_1$ from $I$, and  we obtain two closed subintervals {$L_1$ and $R_1$}. In this case we call $I$ the \emph{father interval} of $U_1$, and  call {$L_1, R_1$} the \emph{generating intervals} of $U_1$.  In the second step, we remove $U_2$. Without loss of generality we may assume $U_2\subset L_1$. Then we obtain two subintervals {$L_2$ and $R_2$ from $L_1$}. This procedure can be continued. Suppose in the $n$-th step we remove $U_n$ from some closed interval {$L_j$ for some $1\le j\le n-1$}, and we obtain two subintervals $L_{n}, R_n$. So {$L_j$} is the father interval of $U_n$, and $L_n, R_n$ are the two generating subintervals of $U_n$. Continuing this procedure indefinitely we  obtain the Cantor set $E$. Then the thickness of $E$ introduced by Newhouse \cite{Newhouse-1970} is defined as
\begin{equation}\label{eq:thickness-newhouse}
\tau_\mathcal N(E):=\sup\inf_{n\in\N}\set{\frac{|L_n|}{|U_n|}, \frac{|R_n|}{|U_n|}},
\end{equation}
where $L_n$ and $R_n$ are the generating intervals of $U_n$ in the procedure, and $|U|$ denotes the length of a subinterval $U\subset\R$. Here the supermum is taken over all permutations of $\set{U_i}_{i=1}^\f$. It is worth to mention that the supermum in (\ref{eq:thickness-newhouse}) can be achieved by ordering the lengths of the open intervals $\set{U_i}_{i=1}^\f$ in a decreasing order (cf.~\cite{Astels_2000}).

Let $x\in(0,1)$, and  let $(x_i)=\Phi_x(1/m)$ be the greedy $m$-adic  expansion of $x$. Then there exist infinitely many $i\ge 1$ such that $x_i<m-1$. Denote by $(n_j)\subset\N$ the set of  all indices $i\ge 1$ such that  $x_i<m-1$. Then   for any $j\ge 1$ we have $x_{n_j}<m-1$ and  $x_{n_j+1}\ldots x_{n_{j+1}-1}=(m-1)^{n_{j+1}-n_j-1}$. For $j\ge 1$ and $b\in\set{x_{n_j}+1, x_{n_j}+2,\ldots, m-1}$ let $p_{j,b}, q_{j,b}\in(0,1/m)$ {be defined by}
\begin{equation}\label{eq:12-1}
\Phi_x(p_{j,b})=x_1\ldots x_{n_j-1}\;b\;(m-1)^\f\quad\textrm{and}\quad \Phi_x(q_{j,b})=x_1\ldots x_{n_j-1}\;b\; 0^\f.
\end{equation}
Then by Lemma \ref{lem:lexicographical-order} it follows that
  \begin{equation}\label{eq:intervalsequence}
  \begin{split}
  0&<p_{j,m-1}<q_{j,m-1}<p_{j,m-2}<q_{j,m-2}<\cdots<p_{j,x_{n_j}+1}<q_{j,x_{n_j}+1}\\
  &<p_{j+1,m-1}<q_{j+1,m-1}<\cdots<p_{j+1,x_{n_{j+1}}+1}<q_{j+1,x_{n_{j+1}}+1}\\
  &<\cdots<1/m
  \end{split}
  \end{equation}
   for all $j\ge 1$, and $p_{j,m-1}\nearrow 1/m$ as $j\to\f$. Note that the intervals $I_{j,b}=[p_{j,b},q_{j,b}]$ with $j\ge 1$ and $b\in\set{x_{n_j}+1,\ldots, m-1}$ are pairwise disjoint.  Then we can rename these intervals  in an increasing order:
      \[I_1=[p_1,q_1], \quad I_2=[p_2, q_2],\quad \cdots, \quad I_k=[p_k,q_k],\quad {I_{k+1}}=[p_{k+1}, q_{k+1}]\quad\cdots\]
   such that $q_k<p_{k+1}$ for all $k\ge 1$. So for each interval $I_k$ there exist a unique $j\in\N$ and a unique $b\in\set{x_{n_j}+1,\ldots, m-1}$ such that $I_k=I_{j,b}$. Since we have only finitely many choices for the index $b$, it follows that $k\to\f$ is equivalent to $j\to\f$.

  Let $k\ge 1$. Then   $I_k=I_{j,b}$ for some $j\ge 1$ and $b\in\set{x_{n_j}+1, \ldots, m-1}$. Now for a word $\w\in\set{0,1,\ldots, m-1}^*$ we define the basic interval associated to $\w$ by
\[
I_k(\w)=[p_{k}(\w), q_{k}(\w)],
\]
where
\begin{equation}\label{eq:thick-1}
\Phi_x(p_{k}(\w))=x_1\ldots x_{n_j-1}b\; \w(m-1)^\f,\quad \Phi_x(q_{k}(\w))=x_1\ldots x_{n_j-1}b \;\w 0^\f.
\end{equation}
Then for the empty word $\epsilon$ we have $I_{k}(\epsilon)=[p_{k}, q_{k}]$.  Observe that  for any two words $\u$ and $\v$, if $\u$ is a prefix of $\v$, then $I_{k}(\v)\subset I_{k}(\u)$. So,  for any $\w\in\set{0,1,\ldots,m-1}^*$  we have $I_{k}(\w)\supset\bigcup_{d=0}^{m-1}I_{k}(\w d)$ with the union pairwise disjoint. Furthermore, the left endpoint of $I_{k}(\w(m-1))$ coincides with the left endpoint of $I_{k}(\w)$, and the right endpoint of $I_{k}(\w0)$ coincides with the right endpoint of $I_{k}(\w)$.  In other words, these intervals $I_{k}(\w)$, $\w\in\set{0,1,\ldots, m-1}^*$ have a tree structure.  Therefore,
\begin{equation}\label{eq:thick-2}
E_{k}(x):=\bigcap_{n=0}^\f\bigcup_{\w\in\set{0,1,\ldots,m-1}^n}I_{k}(\w)
\end{equation}
is a Cantor subset of $\La(x)$. Here  the inclusion $E_{k}(x)\subset \La(x)$ follows by our construction of $E_{k}(x)$ that each $\la\in E_{k}(x)$ corresponds to a  unique coding $(d_i)\in\set{0,1,\ldots,m-1}^\N$ such that $\Phi_x(\la)=(d_i)$, which implies that $\la\in\La(x)$. Since $E_{k}(x)\subset \La(x)$ for all $k\ge 1$, we then construct a sequence of Cantor subsets of $\La(x)$. Furthermore,  by (\ref{eq:intervalsequence}) it follows that these Cantor subsets $E_{k}(x), k\ge 1$ are pairwise disjoint, and $d_H( E_{k}(x), \set{1/m})\to 0$ as $k\to\f$, where $d_H$ is the Hausdorff metric.

 \begin{figure}[h!]
\begin{center}
\begin{tikzpicture}[
    scale=10,
    axis/.style={very thick, ->},
    important line/.style={thick},
    dashed line/.style={dashed, thin},
    pile/.style={thick, ->, >=stealth', shorten <=2pt, shorten
    >=2pt},
    every node/.style={color=black}
    ]

    \draw[important line] (0, 0.7)--(0.4, 0.7);
     \node[] at(0.2, 0.75){$I_{k}(\w^+)$};  \node[] at(0, 0.65){$\la_1=p_{k}(\w^+)$}; \node[] at(0.4, 0.65){$\la_2=q_{k}(\w^+)$};
     \node[] at(0.575, 0.75){$G_{k}(\w)$};
    \draw[important line] (0.7, 0.7)--(1.2, 0.7);
          \node[] at(0.95, 0.75){$I_{k}(\w)$};

 \node[] at(0.7, 0.65){$\la_3=p_{k}(\w)$}; \node[] at(1.2, 0.65){$\la_4=q_{k}(\w)$};
\end{tikzpicture}
\end{center}
\caption{The geometrical structure of the basic intervals $I_{k}(\w^+)=[p_{k}(\w^+), q_{k}(\w^+)]=[\la_1,\la_2]$ and $I_{k}(\w)=[p_{k}(\w), q_{k}(\w)]=[\la_3, \la_4]$, and the gap $G_k(\w)=(q_k(\w^+), p_k(\w))=(\la_2, \la_3)$. }\label{fig:1}
\end{figure}

Let $k\ge 1$. Recall that
for a word $\w=w_1\ldots w_n$ with $w_n<m-1$ we write $\w^+=w_1\ldots w_{n-1}(w_n+1)$.
For two neighboring basic intervals $I_{k}(\w^+)=[p_{k}(\w^+), q_{k}(\w^+)]$ and $I_{k}(\w)=[p_{k}(\w), q_{k}(\w)]$ of the same level,  we call the open interval $G_{k}(\w):=(q_{k}(\w^+), p_{k}(\w))$ the \emph{gap} between them (see Figure \ref{fig:1}).
{Now we write}
\begin{equation}\label{eq:thickness-j}
\begin{split}
\tau(E_{k}(x))&:=\inf_{n\ge 1}\min_{\w,\w^+\in\set{0,1,\ldots, m-1}^n}\set{\frac{|I_{k}(\w^+)|}{|G_{k}(\w)|},\frac{|I_{k}(\w)|}{|G_{k}(\w)|}} \\
&=\inf_{n\ge 1}\min_{\w,\w^+\in\set{0,1,\ldots, m-1}^n}\set{\frac{q_{k}(\w^+)-p_{k}(\w^+)}{p_{k}(\w)-q_{k}(\w^+)}, \frac{q_{k}(\w)-p_{k}(\w)}{p_{k}(\w)-q_{k}(\w^+)}},
\end{split}
\end{equation}

{In the next lemma we show that  $\tau_{\mathcal N}(E_{k}(x))\ge \tau(E_{k}(x))$.}
\begin{lemma}\label{lem:thickness-two-defs}
Let $x\in(0,1)$. Then for  any $k\ge 1$   we have $\tau_{\mathcal N}(E_{k}(x))\ge \tau(E_{k}(x))$.
\end{lemma}
\begin{proof}
Let $x\in(0,1)$  and $k\ge 1$. We remove the open intervals (gaps) $G_{k}(\w)$ with $\w\in\set{0,1,\ldots, m-1}^*$ from $I_k(\epsilon)=[p_k, q_k]$ in the following way. First we remove from $I_k(\epsilon)$ the open intervals $G_{k}(0)$, and then the open interval $G_{k}(1)$, and  so on,  and in the $(m-1)$-th step we remove the open interval $G_{k}(m-2)$. Then by the definitions of  thickness  defined in (\ref{eq:thickness-newhouse}) and (\ref{eq:thickness-j}) respectively  it follows that
\[
  \min_{0\le d\le m-2}\set{\frac{|L_{k}(d)|}{|G_{k}(d)|}, \frac{|R_{k}(d)|}{|G_{k}(d)|}}
\ge \min_{0\le d\le m-2}\set{\frac{|I_{k}(d+1)|}{|G_{k}(d)|}, \frac{|I_{k}(d)|}{|G_{k}(d)|}},\]
where $L_{k}(d)$ and $R_{k}(d)$ are the generating intervals of $G_{k}(d)$ in the construction of $E_{k}(x)$ by Newhouse. After removing the $(m-1)$ open intervals $G_k(0), G_k(1), \ldots, G_k(m-2)$ we obtain $m$ basic intervals $I_k(0), I_k(1),\ldots, I_k(m-1)$.  Next for each basic interval $I_k(d)$ we remove from $I_k(d)$ the open  intervals $G_{k}(d0), G_k(d1),\ldots, G_k(d(m-2))$ successively, and get $m$ basic subintervals $I_k(d0), I_k(d1), \ldots, I_{k}(d(m-1))$.  Proceeding this argument, and suppose we are considering the basic interval  $I_{k}(\w)$ for some $\w$ of length $n$. Then we successively remove from $I_k(\w)$ the open intervals $G_{k}(\w0), G_{k}(\w1),\ldots, G_{k}(\w(m-2))$, and obtain $m$ basic subintervals $I_k(\w 0), I_k(\w 1), \ldots, I_k(\w(m-1))$. By the same argument as above it is easy to see that
\[
\min_{0\le d\le m-2}\set{\frac{|L_{k}(\w d)|}{|G_{k}(\w d)|}, \frac{|R_{k}(\w d)|}{|G_{k}(\w d)|}}
\ge \min_{0\le d\le m-2}\set{\frac{|I_{k}(\w (d+1))|}{|G_{k}(\w d)|}, \frac{|I_{k}(\w d)|}{|G_{k}(\w d)|}},
\]
where $L_{k}(\w d)$ and $R_{k}(\w d)$ are the generating intervals of $G_{k}(\w d)$. By induction  it follows that
\[
\min_{\w,\w^+\in\set{0,1,\ldots, m-1}^n}\set{\frac{|L_{k}(\w)|}{|G_{k}(\w)|}, \frac{|R_{k}(\w)|}{|G_{k}(\w)|}}
\ge \min_{\w,\w^+\in\set{0,1,\ldots, m-1}^n}\set{\frac{|I_{k}(\w^+)|}{|G_{k}(\w)|}, \frac{|I_{k}(\w)|}{|G_{k}(\w)|}}
\]
for all $n\ge 1$. Therefore, by (\ref{eq:thickness-newhouse}) and (\ref{eq:thickness-j}) we conclude that
\begin{align*}
\tau_{\mathcal N}(E_{k}(x))&\ge \inf_{n\ge 1}\min_{\w,\w^+\in\set{0,1,\ldots, m-1}^n}\set{\frac{|L_{k}(\w)|}{|G_{k}(\w)|}, \frac{|R_{k}(\w)|}{|G_{k}(\w)|}}\\
&\ge \inf_{n\ge 1}\min_{\w,\w^+\in\set{0,1,\ldots, m-1}^n}\set{\frac{|I_{k}(\w^+)|}{|G_{k}(\w)|}, \frac{|I_{k}(\w)|}{|G_{k}(\w)|}}=\tau(E_{k}(x)).
\end{align*}
\end{proof}

In the following we   show that  the thickness of $E_k(x)$ goes to infinitey as $k\to\f$.
\begin{proposition}\label{prop:thickness-La(x)}
Let $x\in(0,1)$. Then $\tau(E_k(x))\to +\f$ as $k\to\f$.
\end{proposition}
\begin{proof}
Let $x\in(0,1)$, {and let $(x_i)=\Phi_x(1/m)$ be the greedy $m$-adic expansion of $x$}. Then there exists a smallest integer $\ell\ge 1$ such that $x_\ell>0$.
Take $k$ large enough, so   there exist  $j\ge 1, b\in\set{x_{n_j}+1, \ldots, m-1}$ such that  $I_k=I_{j,b}$ and  $n_j>\ell$.
For $n\ge 1$ let $\w\in\set{0,1,\ldots, m-1}^n$ such that $\w^+\in\set{0,1,\ldots, m-1}^n$. In view of (\ref{eq:thickness-j}) we need to estimate the lower bounds for the two quotients
  \[
  \frac{q_{k}(\w^+)-p_{k}(\w^+)}{p_{k}(\w)-q_{k}(\w^+)}, \quad \frac{q_{k}(\w)-p_{k}(\w)}{p_{k}(\w)-q_{k}(\w^+)}.
  \]
  For simplicity we write $\la_1=p_{k}(\w^+), \la_2=q_{k}(\w^+), \la_3=p_{k}(\w)$ and $\la_4=q_{k}(\w)$. Then $\la_1<\la_2<\la_3<\la_4$ (see Figure \ref{fig:1}).
  So we need to  estimate the lower bounds for $\la_2-\la_1$ and $\la_4-\la_3$, and the upper bounds for $\la_3-\la_2$.

  {\bf A lower bound for $\la_2-\la_1$.}  Note by (\ref{eq:thick-1}) that
  \[
  \pi_{\la_1}(x_1\ldots x_{n_j-1}b \w^+ (m-1)^\f)=x=\pi_{\la_2}(x_1\ldots x_{n_j-1}b \w^+ 0^\f).
  \]
    Rearranging the above equation we obtain
  \begin{align*}
   \pi_{\la_2} (0^{n_j+n}(m-1)^\f)&=\pi_{\la_2}(x_1\ldots x_{n_j-1}b \w^+(m-1)^\f)-\pi_{\la_1}(x_1\ldots x_{n_j-1}b \w^+ (m-1)^\f)\\
   &\le\pi_{\la_2}((m-1)^\f)-\pi_{\la_1}((m-1)^\f)\\
   &=\frac{m-1}{1-\la_2}\la_2-\frac{m-1}{1-\la_1}\la_1\\
   &=\frac{m-1}{(1-\la_1)(1-\la_2)}(\la_2-\la_1).
  \end{align*}
  This implies that
  \begin{equation}
    \label{eq:09-1}
    \la_2-\la_1\ge \frac{(1-\la_1)(1-\la_2)}{m-1}\pi_{\la_2}(0^{n_j+n}(m-1)^\f)=(1-\la_1)\la_2^{n_j+n+1}.
  \end{equation}

  {\bf A lower bound for $\la_4-\la_3$.}
  Note by (\ref{eq:thick-1}) that
   \[
  \pi_{\la_3}(x_1\ldots x_{n_j-1}b \w (m-1)^\f)=x=\pi_{\la_4}(x_1\ldots x_{n_j-1}b \w  0^\f).
  \]
  Then by the same argument as in the estimate for $\la_2-\la_1$ one can verify that
  \begin{equation}
    \label{eq:09-2}
    \la_4-\la_3\ge (1-\la_3)\la_4^{n_j+n+1}.
  \end{equation}

  {\bf An upper bound for $\la_3-\la_2$.}
  Note by (\ref{eq:thick-1}) that
  \[\pi_{\la_2}(x_1\ldots x_{n_j-1}b \w^+ 0^\f)=x=\pi_{\la_3}(x_1\ldots x_{n_j-1}b \w(m-1)^\f).\] This implies that
  \begin{align*}
    &\pi_{\la_2}(0^{n_j+n-1}10^\f)-\pi_{\la_3}(0^{n_j+n}(m-1)^\f)\\
    =&\pi_{\la_3}(x_1\ldots x_{n_j-1}b \w0^\f)-\pi_{\la_2}(x_1\ldots x_{n_j-1}b \w0^\f)\\
    \ge& \la_3^\ell-\la_2^\ell\ge\la_2^{\ell-1}(\la_3-\la_2),
  \end{align*}
  where the first inequality follows by the definition of $\ell$ that $\ell<n_j$ and $x_\ell>0$. So, by using $\la_3>\la_2$ it follows that
  \begin{equation}\label{eq:09-3}
  \begin{split}
    \la_3-\la_2&\le \la_2^{1-\ell}\left(\la_2^{n_j+n}-\frac{m-1}{1-\la_3}\la_3^{n_j+n+1}\right)\\
    &\le \la_2^{n_j-\ell+n+1}\left(1-\frac{m-1}{1-\la_3}\la_3\right)\\
    &=\la_2^{n_j-\ell+n+1}\frac{m}{1-\la_3}\left(\frac{1}{m}-\la_3\right).
  \end{split}
  \end{equation}
  From this we still need to estimate $1/m-\la_3$. Note that $\Phi_x(1/m)=(x_i)$. Then
  \[
  \pi_{\la_3}(x_1\ldots x_{n_j-1}b \w (m-1)^\f)=x=\pi_{1/m}(x_1x_2\ldots).
  \]
  This yields that
  \begin{align*}
    &\pi_{\la_3}(0^{n_j-1}(b-x_{n_j}) \w(m-1)^\f)-\pi_{1/m}(0^{n_j}x_{n_j+1} x_{n_j+2}\ldots)\\
    =&\pi_{1/m}(x_1\ldots  x_{n_j}0^\f)-\pi_{\la_3}(x_1\ldots x_{n_j}0^\f)\\
    \ge& (1/m)^\ell-\la_3^\ell\ge \la_3^{\ell-1}(1/m-\la_3).
  \end{align*}
  So, by using $b\in\set{x_{n_j}+1,\ldots, m-1}$ it follows that
  \begin{equation}
    \label{eq:09-4}
    \begin{split}
    1/m-\la_3&\le \la_3^{1-\ell}\pi_{\la_3}(0^{n_j-1}(b-x_{n_j})\w(m-1)^\f)\\
    &\le\la_3^{1-\ell}\pi_{\la_3}(0^{n_j-1}(m-1)^\f)\\
    &<\la_3^{1-\ell}\la_3^{n_j-1}=\la_3^{n_j-\ell}.
    \end{split}
  \end{equation}
Substituting (\ref{eq:09-4}) into (\ref{eq:09-3}) we obtain an upper bound for $\la_3-\la_2$:
  \begin{equation}
    \label{eq:09-5}
    \la_3-\la_2\le   \frac{m\la_3^{n_j-\ell}}{1-\la_3}\la_2^{n_j-\ell+n+1}.
  \end{equation}

  Therefore, by (\ref{eq:09-1}) and (\ref{eq:09-5}) it follows that
  \begin{equation}
    \label{eq:09-6}
     \frac{q_{k}(\w^+)-p_{k}(\w^+)}{p_{k}(\w)-q_{k}(\w^+)}=\frac{\la_2-\la_1}{\la_3-\la_2}\ge \frac{(1-\la_1)(1-\la_3)\la_2^\ell}{m\la_3^{n_j-\ell}}\to +\f \quad\textrm{as }k\to\f.
  \end{equation}
  Here we emphasize that $\ell\in\N$ depends only on $x$, and $k\to\f$ implies $j\to\f$. Note that $\la_4>\la_2$. Then   by (\ref{eq:09-2}) and (\ref{eq:09-5}) we obtain that
  \begin{equation}
    \label{eq:09-7}
    \frac{q_{k}(\w)-p_{k}(\w)}{p_{k}(\w)-q_{k}(\w^+)}=\frac{\la_4-\la_3}{\la_3-\la_2}\ge  \frac{(1-\la_3)^2\la_4^\ell}{m\la_3^{n_j-\ell}}\to +\f \quad\textrm{as }k\to\f.
  \end{equation}
  Hence, {according to  (\ref{eq:thickness-j})}, and by (\ref{eq:09-6}) and (\ref{eq:09-7}) we conclude that $\tau(E_{k}(x))\to +\f$ as $k\to\f$. This completes the proof.
\end{proof}
\begin{remark}\label{rem:1}
  \begin{enumerate}

  \item Note  by \cite{Palis_Takens_1993} that for any Cantor set $E$ we have
  \begin{equation}\label{eq:dim-thickness}
 \dim_H E\ge\frac{\log 2}{\log (2+\frac{1}{\tau_{\mathcal N}(E)})}.
  \end{equation}
  So, by Lemma \ref{lem:thickness-two-defs}, Proposition \ref{prop:thickness-La(x)} and (\ref{eq:dim-thickness}) it follows that  for any $x\in(0,1)$,
  \[\dim_H\La(x)\ge\dim_H E_k(x)\ge\frac{\log 2}{\log(2+\frac{1}{\tau_{\mathcal N}(E_k(x))})}\ge \frac{\log 2}{\log(2+\frac{1}{\tau(E_k(x))})}\to 1 \quad\textrm{as }k\to\f.\]   This provides an alternative  proof of $\dim_H\La(x)=1$.

\item {It is well-known that for a  Cantor sets $E\subset \R$, if $\tau_{\mathcal N}(E)\ge 1$, then the algebraic sum $E+E:=\set{a+b: a, b\in E}$ contains an interval (cf.~\cite{Newhouse-1970}).  So, Proposition \ref{prop:thickness-La(x)} suggests that {the} algebraic sum $\La(x)+\La(x)$
contains an interval for any $x\in(0,1)$.}

\item Recall that a set $E\subset\R$ is said to contain \emph{arbitrarily long arithmetic progressions} if  {for any $n\in\N$} there exist $a, b\in\R$ such that $\set{a+b,a+2b,a+3b,\ldots, a+n b}\subset E$. Recently, Yavicoli \cite[Remark of Theorem 4]{Yavicoli-2020} proved that for any $n\in\N$ there exists a (finite) constant $\tau_n$ such that  if the thickness of $E$ is larger than $\tau_n$, then $E$ contains   arithmetic progressions of length $n$. This, together with Proposition \ref{prop:thickness-La(x)}, implies that $\La(x)$ contains arbitrarily long arithmetic progressions {for any $x\in(0,1)$}.

\end{enumerate}

\end{remark}

\subsection{Large intersections of  $\La(x)$ and $\La(y)$}
In order to prove Theorem \ref{th:intersection-La(x)-La(y)} we recall the following results from  Kraft \cite[Theorem 1.1]{Kraft-1992} and Hunt et al.~\cite[Theorem 1]{Hunt-Kan-Yorke-1992}. Two Cantor sets $F_1$ and $F_2$ are called \emph{interleaved} if \[F_1\cap conv(F_2)\ne\emptyset\quad \textrm{and}\quad conv(F_1)\cap F_2\ne\emptyset.\]
\begin{lemma}\label{lem:thickness-intersection}
  Let $x, y\in(0,1)$. If there exist $i, j\in\N$  such that $E_{i}(x)$ and $E_{j}(y)$ are interleaved with
  $\tau_{i,j}(x,y):=\min\set{\tau_{\mathcal N}(E_{i}(x)), \tau_{\mathcal N}(E_{j}(y))}>1+\sqrt{2}$,
  then $E_{i}(x)\cap E_{j}(y)$ contains a Cantor subset of thickness at least  of order $\sqrt{\tau_{i,j}(x,y)}$.
\end{lemma}
By Lemma \ref{lem:thickness-two-defs} and the proof of  Proposition \ref{prop:thickness-La(x)} it follows that for any $x\in(0,1)$ we have $\tau_{\mathcal N}(E_i(x))\to \f$ as $i\to\f$. So, if we can show that there exist infinitely many pairs $(i, j)\in\N\times\N$ such that $E_i(x)$ and $E_j(y)$ are interleaved Cantor sets, then Theorem \ref{th:intersection-La(x)-La(y)} can be deduced from Lemma \ref{lem:thickness-intersection} and (\ref{eq:dim-thickness}).
In the following we will show that  there exist infinitely many pairs $(i,j)$ such that $E_{i}(x)$ and $E_{j}(y)$ are interleaved.

\begin{figure}[h!]
\begin{center}
\begin{tikzpicture}[
    scale=10,
    axis/.style={very thick, ->},
    important line/.style={thick},
    dashed line/.style={dashed, thin},
    pile/.style={thick, ->, >=stealth', shorten <=2pt, shorten
    >=2pt},
    every node/.style={color=black}
    ]

    \draw[important line] (0, 0.7)--(0.4, 0.7);
     \node[] at(0.2, 0.75){$I_{k}$};  \node[] at(0, 0.65){$p_{k}$}; \node[] at(0.4, 0.65){$q_{k}$};
    \draw[important line] (0.7, 0.7)--(1.0, 0.7);
          \node[] at(0.85, 0.75){$I_{k+1}$};
 \node[] at(0.7, 0.65){$p_{k+1}$}; \node[] at(1.0, 0.65){$q_{k+1}$};

          \draw[dashed line] (1.1,0.7)--(1.3,0.7);
  \node[] at(1.3,0.75){$\frac{1}{m}$};

\end{tikzpicture}
\end{center}
\caption{The geometrical structure of the basic intervals $I_k=[p_{k},q_{k}]=conv(E_k(x))$  and $I_{k+1}=[p_{k+1},q_{k+1}]=conv(E_{k+1}(x))$. They converge to $\set{1/m}$ in the Hausdorff metric $d_H$ as $k\to\f$.}\label{fig:2}
\end{figure}
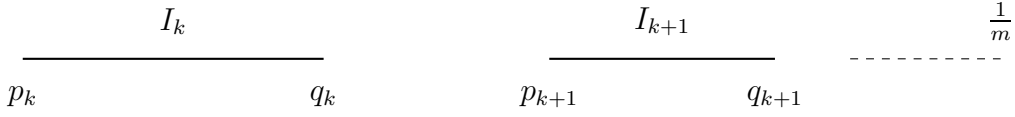

For $x\in(0,1)$ we recall from (\ref{eq:thick-2}) that $E_{k}(x)$ is a Cantor subset of $\La(x)$ for any $k\ge 1$. For each basic interval   $I_k=conv(E_{k}(x))=[p_{k}, q_{k}]$ there exist a unique $j\ge 1$ and $b\in\set{x_{n_j}+1,\ldots,m-1}$ such that $I_k=I_{j,b}$. Then the endpoints $p_k, q_k$ are determined by
\[
\Phi_x(p_{k})=x_1\ldots x_{n_j-1}b\;(m-1)^\f\quad\textrm{and}\quad \Phi_x(q_{k})=x_1\ldots x_{n_j-1}b\; 0^\f.
\] First we define the \emph{thickness} of  the interval sequence  $\set{I_{k}:k\ge 1}$ by
\[
\theta(x):=\liminf_{k\to\f}\theta_{k}(x),
\]
where
\begin{equation}\label{eq:theta-ib}
\theta_{k}(x):=
\min\set{\frac{|I_{k}|}{|G_{k}|}, \frac{|I_{k+1}|}{|G_{k}|}}=\min\set{\frac{q_k-p_k}{p_{k+1}-q_k}, \frac{q_{k+1}-p_{k+1}}{p_{k+1}-q_k}}.
\end{equation}
Here $G_{k}=(q_k, p_{k+1})$ is the gap between   $I_{k}$ and $I_{k+1}$ (see Figure \ref{fig:2}).

In the following lemma we show that  the thickness of the interval sequence $\set{I_k: k\ge 1}$ is infinity for any $x\in(0,1)$.
\begin{lemma}
  \label{lem:thickness-interval-sequences}
 For any $x\in(0,1)$ we have $\theta(x)=+\f$.
\end{lemma}
\begin{proof}
  The idea to prove this lemma is similar to that for Proposition \ref{prop:thickness-La(x)}. By the definition of $\theta(x)$ we need to estimate the lower bounds $\theta_{k}(x)$ for $k$ large enough. Take $k\in\N$ sufficiently large. Then  there exist a unique $j\in\N$ and a unique $b\in\set{x_{n_j}+1,\ldots, m-1}$ such that $I_k=I_{j,b}$. In view of (\ref{eq:theta-ib}) we consider the following two cases: (I) $b\in\set{x_{n_j}+2,\ldots, m-1}$; (II) $b=x_{n_j}+1$.

  Case (I). $b\in\set{x_{n_j}+2,\ldots, m-1}$. Then $I_k=I_{j,b}$ and $I_{k+1}=I_{j,b-1}$.  We need to estimate the lower bounds of $|I_k|=q_{k}-p_{k}$ and $|I_{k+1}|=q_{k+1}-p_{k+1}$, and the upper bounds of $|G_k|=p_{k+1}-q_{k}$.
  Observe that
  \[
  \pi_{q_{k}}(x_1\ldots x_{n_j-1}b 0^\f)=x=\pi_{p_k}(x_1\ldots x_{n_j-1}b(m-1)^\f).
  \]
  This implies that
  \begin{align*}
    \pi_{q_{k}}(0^{n_j}(m-1)^\f)&=\pi_{q_{k}}(x_1\ldots x_{n_j-1}b (m-1)^\f)-\pi_{p_{k}}(x_1\ldots x_{n_j-1}b (m-1)^\f)\\
    &\le \frac{m-1}{1-q_{k}}q_{k}-\frac{m-1}{1-p_{k}}p_{k}\\
    &=\frac{m-1}{(1-p_{k})(1-q_{k})}(q_{k}-p_{k}).
  \end{align*}
  Whence,
  \begin{equation}
    \label{eq:13-1}
    q_{k}-p_{k}\ge \frac{(1-p_{k})(1-q_{k})}{m-1}\times\frac{m-1}{1-q_{k}}q_{k}^{n_j+1}=(1-p_{k})q_{k}^{n_j+1}.
  \end{equation}
Similarly, one can prove that
 \begin{equation}
    \label{eq:13-2}
    q_{k+1}-p_{k+1}\ge  (1-p_{k+1})q_{k+1}^{n_{j}+1}.
  \end{equation}

  Now we turn to the upper bounds of  $p_{k+1}-q_{k}$. Note that
  \[
  \pi_{q_{k}}(x_1\ldots x_{n_j-1}b 0^\f)=x=\pi_{p_{k+1}}(x_1\ldots x_{n_j-1}(b-1) (m-1)^\f).
  \]
  Then
  \begin{align*}
   & \pi_{q_{k}}(0^{n_j-1}1 0^\f)-\pi_{p_{k+1}}(0^{n_j}(m-1)^\f) \\
       =&\pi_{p_{k+1}}(x_1\ldots x_{n_j-1}(b-1) 0^\f)-\pi_{q_{k}}(x_1\ldots x_{n_j-1}(b-1) 0^\f)\\
    \ge& p_{k+1}^\ell- q_{k}^{\ell}\ge q_{k}^{\ell-1}(p_{k+1}-q_{k}),
  \end{align*}
  where $\ell<n_j$ is the smallest integer such that $x_\ell>0$. Therefore,
  \begin{equation}\label{eq:13-3}
  \begin{split}
    p_{k+1}-q_{k}&\le q_{k}^{1-\ell}\left(q_{k}^{n_j}-\frac{m-1}{1-p_{k+1}}p_{k+1}^{n_j+1}\right)\\
    &\le q_{k}^{n_j-\ell+1}\left(1-\frac{m-1}{1-p_{k+1}}p_{k+1}\right)\\
    &=q_{k}^{n_j-\ell+1}\frac{m}{1-p_{k+1}}\left(\frac{1}{m}-p_{k+1}\right).
  \end{split}
  \end{equation}
 So we still need to estimate $1/m-p_{k+1}$. Observe that
 \[
 \pi_{p_{k+1}}(x_1\ldots x_{n_{j}-1}(b-1) (m-1)^\f)=x=\pi_{1/m}(x_1\ldots x_{n_{j}}x_{n_{j}+1}\ldots).
 \]
 Then
 \begin{align*}
   \pi_{p_{k+1}}(0^{n_{j}-1} (b-1-x_{n_j})(m-1)^\f) &\ge \pi_{1/m}(x_1\ldots x_{n_{j}}0^\f)-\pi_{p_{k+1}}(x_1\ldots x_{n_{j}}0^\f)\\
   &\ge (1/m)^\ell-p_{k+1}^\ell\ge p_{k+1}^{\ell-1}(1/m-p_{k+1}),
 \end{align*}
 which implies that
 \begin{equation}
   \label{eq:13-4}
   1/m-p_{k+1}\le p_{k+1}^{1-\ell}  p_{k+1}^{n_{j}-1}=p_{k+1}^{n_{j}-\ell}.
 \end{equation}
 Substituting (\ref{eq:13-4}) into (\ref{eq:13-3}) we obtain an upper bound for $p_{k+1}-q_k$:
 \begin{equation}
   \label{eq:13-5}
   p_{k+1}-q_{k}\le q_{k}^{n_j-\ell+1}\frac{m}{1-p_{k+1}} p_{k+1}^{n_{j}-\ell}.
 \end{equation}

 Hence, by (\ref{eq:13-1}), (\ref{eq:13-2}) and (\ref{eq:13-5}) it follows that
 \begin{align*}
  \theta_{k}(x)=&\min\set{\frac{q_{k}-p_{k}}{p_{k+1}-q_{k}}, \frac{q_{k+1}-p_{k+1}}{p_{k+1}-q_{k}}}\\
   \ge& \min\set{\frac{(1-p_{k})(1-p_{k+1})q_{k}^\ell}{m p_{k+1}^{n_{j}-\ell}}, \frac{(1-p_{k+1})^2 q_{k+1}^\ell}{m p_{k+1}^{n_j-\ell}}}\to\;+\f
 \end{align*}
  as $k\to\f$.

  Case (II). $b=x_{n_j}+1$. Then $I_k=I_{j,b}=[p_k,q_k]$ and $I_{k+1}=I_{j+1,m-1}=[p_{k+1}, q_{k+1}]$. In this case we have
  \begin{align*}
  \Phi_x(p_k)&=x_1\ldots x_{n_j}^+(m-1)^\f,\quad \Phi_x(q_k)=x_1\ldots x_{n_j}^+ 0^\f;\\
  \Phi_x(p_{k+1})&=x_1\ldots x_{n_{j+1}-1}(m-1)^\f=x_1\ldots x_{n_j}(m-1)^\f,\\
   \Phi_x(q_{k+1})&=x_1\ldots x_{n_{j+1}-1}(m-1)0^\f=x_1\ldots x_{n_j}(m-1)^{n_{j+1}-n_j} 0^\f.
  \end{align*}
  Here we emphasize that by the definition of $(n_j)$ it follows that $x_i=m-1$ for any $n_j<i<n_{j+1}$.
    By the same argument as in  Case I one can prove that
  \begin{equation}\label{eq:kong1}
  q_{k}-p_{k}\ge (1-p_{k})q_{k}^{n_j+1},\quad
  q_{k+1}-p_{k+1}\ge (1-p_{k+1})q_{k+1}^{n_{j+1}+1},
  \end{equation}
  and
  \begin{equation}\label{eq:kong2}
  p_{k+1}-q_{k}\le q_{k}^{n_j-\ell+1}\frac{m}{1-p_{k+1}}p_{k+1}^{n_{j+1}-\ell}.
  \end{equation}
  Therefore,  by (\ref{eq:kong1}) and (\ref{eq:kong2}) it follows that
  \[
  \theta_{k}(x)=\min\set{\frac{q_{k}-p_{k}}{p_{k+1}-q_{k}}, \frac{q_{k+1}-p_{k+1}}{p_{k+1}-q_{k}}}\to+\f
  \]
  as $k\to\f$.

  Hence, by Cases (I) and (II) we conclude that $\theta(x)=+\f$.
\end{proof}

In order to prove that $E_{i}(x)$ and $E_{j}(y)$ are interleaved for infinitely many pairs $(i,j)$ we still need the following estimate.
\begin{lemma}
  \label{lem:interval-comparison}
  Let $x\in(0,1)$. Then there exist   a constant $\eta(x)>1$ and a large integer $N$ such that for any $k\ge N$ we have
  \[
  \frac{|I_{k+1}|}{|I_{k}|}=\frac{q_{k+1}-p_{k+1}}{q_{k}-p_{k}}<\eta(x).
    \]
\end{lemma}

\begin{proof}
  Suppose $I_k=I_{j,b}{=[p_k, q_k]}$ for some $j\ge 1$ and $b\in\set{x_{n_j}+1,\ldots, m-1}$. We split the proof into the following two cases: (I) $b\in\set{x_{n_j}+2,\ldots, m-1}$; (II) $b=x_{n_j}+1$.

Case (I). $b\in\set{x_{n_j}+2,\ldots, m-1}$. Then $I_{k+1}=I_{j,b-1}=[p_{k+1}, q_{k+1}]$.  Note by (\ref{eq:13-1})  that $q_k-p_k\ge (1-p_k)q_k^{n_j+1}$. In the following it suffices to estimate the upper bounds of $q_{k+1}-p_{k+1}$.
Observe that
\[
\pi_{q_{k+1}}(x_1\ldots x_{n_{j}-1}(b-1) 0^\f)=x=\pi_{p_{k+1}}(x_1\ldots x_{n_{j}-1}(b-1)(m-1)^\f).
\]
Then
\begin{align*}
  \pi_{p_{k+1}}(0^{n_{j}}(m-1)^\f)&=\pi_{q_{k+1}}(x_1\ldots x_{n_{j}-1}(b-1) 0^\f)-\pi_{p_{k+1}}(x_1\ldots x_{n_{j}-1}(b-1) 0^\f)\\
  &\ge q_{k+1}^\ell-p_{k+1}^\ell\ge p_{k+1}^{\ell-1}(q_{k+1}-p_{k+1}),
\end{align*}
where $\ell<n_{j}$ is the smallest integer such that $x_\ell>0$. This implies that
\begin{equation}
  \label{eq:15-1}
  q_{k+1}-p_{k+1}\le p_{k+1}^{1-\ell}\frac{m-1}{1-p_{k+1}}p_{k+1}^{n_{j}+1}< p_{k+1}^{n_{j}-\ell+1},
\end{equation}
where the last inequality follows by $p_{k+1}<1/m$. Therefore, by (\ref{eq:13-1}) and (\ref{eq:15-1}) it follows that
\begin{equation}
  \label{eq:15-2}
  \frac{q_{k+1}-p_{k+1}}{q_k-p_k}<\frac{p_{k+1}^{n_{j}-\ell+1}}{(1-p_k)q_k^{n_j+1}}=\frac{1}{(1-p_k)p_{k+1}^{\ell}}\left(\frac{p_{k+1}}{q_k}\right)^{n_{j}+1}.
\end{equation}
Since $p_k\nearrow 1/m$ as $k\to\f$, there exists a large  integer $N_1$ such that $p_{k+1}>\frac{1}{2m}$ for any $k>N_1$.   Then (\ref{eq:15-2}) suggests that
\begin{equation}
  \label{eq:15-3}
  \frac{q_{k+1}-p_{k+1}}{q_k-p_k}<\frac{1}{(1-1/m)(1/2m)^{\ell}}\left(\frac{p_{k+1}}{q_k}\right)^{n_{j}+1}=\frac{m(2m)^{\ell}}{m-1}\left(\frac{p_{k+1}}{q_k}\right)^{n_{j}+1}
\end{equation}
for any $k>N_1$.
Note that $k\to\f$ is equivalent to $j\to\f$. So, to finish the proof it suffices to prove that $\lim_{k\to\f}(p_{k+1}/q_k)^{n_{j}+1}=1$. This follows by (\ref{eq:13-5}) that
\begin{align*}
 {\limsup}_{k\to\f} \left(\frac{p_{k+1}}{q_k}\right)^{n_{j}+1}&={\limsup}_{k\to\f}\left(1+\frac{p_{k+1}-q_k}{q_k}\right)^{n_j+1}\\
 &\le {\limsup}_{k\to\f}\left(1+q_k^{n_j-\ell}\frac{m}{1-p_{k+1}}p_{k+1}^{n_{j}-\ell}\right)^{n_{j}+1}\\
 &=\exp\left[{{\limsup}_{k\to\f}\left((n_{j}+1)q_k^{n_j-\ell}\frac{m}{1-p_{k+1}}p_{k+1}^{n_{j}-\ell}\right)}\right]=\; 1,
\end{align*}
and that $(p_{k+1}/q_k)^{n_{j}+1}\ge 1$ for any $k\ge 1$. So, by (\ref{eq:15-3})  there must exist an integer  $N>N_1$ such that for any $k>N$,
\[
 \frac{q_{k+1}-p_{k+1}}{q_k-p_k}<\frac{(2m)^{\ell+1}}{m-1}
\]
as desired.

Case (II). $b=x_{n_j}+1$. Then $I_{k+1}=I_{j+1,m-1}=[p_{k+1}, q_{k+1}]$. Note by (\ref{eq:kong1}) that  $q_k-p_k\ge (1-p_k)q_k^{n_j+1}$. In the following it suffices to estimate the upper bounds for $q_{k+1}-p_{k+1}$. Observe that
\[
\pi_{q_{k+1}}(x_1\ldots x_{n_j}(m-1)^{n_{j+1}-n_j} 0^\f)=x=\pi_{p_{k+1}}(x_1\ldots x_{n_j}(m-1)^\f).
\]
Then by the same argument as in Case (I) one can show that
\[
q_{k+1}-p_{k+1}\le p_{k+1}^{n_{j+1}-\ell+1}.
\]
This together with (\ref{eq:kong1}) implies that
\begin{equation}\label{eq:30-1}
\frac{q_{k+1}-p_{k+1}}{q_k-p_k}\le\frac{p_{k+1}^{n_{j+1}-\ell+1}}{(1-p_k)q_k^{n_j+1}}=\frac{q_k^{n_{j+1}-n_j}}{(1-p_k)p_{k+1}^\ell}\left(\frac{p_{k+1}}{q_k}\right)^{n_{j+1}+1}.
\end{equation}
Since $p_k\nearrow 1/m$ as $k\to\f$, and $n_{j+1}-n_j\ge 1$, there exists a large integer $N_2$ such that $p_{k+1}>1/(2m)$ for any $k> N_2$. Then (\ref{eq:30-1}) suggests that
\begin{equation}\label{eq:30-2}
\frac{q_{k+1}-p_{k+1}}{q_k-p_k}\le\frac{1}{(1-1/m)(1/2m)^\ell}\left(\frac{p_{k+1}}{q_k}\right)^{n_{j+1}+1}=\frac{m(2m)^\ell}{m-1}\left(\frac{p_{k+1}}{q_k}\right)^{n_{j+1}+1}
\end{equation}
for any $k>N_2$.
Note by (\ref{eq:kong2}) that
\begin{align*}
 {\limsup}_{k\to\f} \left(\frac{p_{k+1}}{q_k}\right)^{n_{j+1}+1}&={\limsup_{k\to\f}}\left(1+\frac{p_{k+1}-q_k}{q_k}\right)^{n_{j+1}+1}\\
 &\le {\limsup}_{k\to\f}\left(1+q_k^{n_j-\ell}\frac{m}{1-p_{k+1}}p_{k+1}^{n_{j+1}-\ell}\right)^{n_{j+1}+1}\\
 &=\exp\left[{{\limsup}_{k\to\f}\left((n_{j+1}+1)p_{k+1}^{n_{j+1}-\ell} q_k^{n_j-\ell}\frac{m}{1-p_{k+1}}\right)}\right]=\; 1.
\end{align*}
So, by (\ref{eq:30-2})  there must exist an integer  $N>N_2$ such that for any $k>N$,
\[
 \frac{q_{k+1}-p_{k+1}}{q_k-p_k}<\frac{(2m)^{\ell+1}}{m-1}.
\]
This completes the proof.
\end{proof}
Based on Lemmas \ref{lem:thickness-interval-sequences} and \ref{lem:interval-comparison} we are ready to show that
 $E_i(x)$ and $E_j(y)$ are interleaved for infinitely many pairs
$(i, j)\in\N\times\N$.
\begin{lemma}
  \label{lem:interleaved-pairs}
 For any $x, y\in(0,1)$ there exists a sequence of pairs $(i_k,j_k)\in\N\times\N$ such that for any $k\ge 1$ we have $i_k<i_{k+1}, j_k<j_{k+1}$, and the two Cantor sets $E_{i_k}(x), E_{j_k}(y)$ are interleaved.
\end{lemma}
\begin{proof}
  Take $x, y\in(0,1)$. Clearly the lemma holds {true} if $x=y$. So in the following  we assume $x\ne y$.    To emphasize the dependence on $x$ we write $I_i^x:=conv(E_i(x))=[p_i^x, q_i^x]$. We also denote by $G_i^x:=(q_i^x, p_{i+1}^x)$ the gap between the two intervals $I_i^x$ and $I_{i+1}^x$. Observe that the closed intervals $\set{I_i^x: i\ge 1}$ are pairwise disjoint, and for each $i\ge 1$ the interval $I_{i+1}^x$ is on the righthand side of $I_i^x$ (see Figure \ref{fig:2}). Furthermore, $p_i^x\nearrow 1/m$   as $i\to\f$. Let $\eta(x)>1$ and $\eta(y)>1$ be the constants defined as in Lemma \ref{lem:interval-comparison}, and let $\eta:=\max\set{\eta(x), \eta(y)}$.  Then by Proposition \ref{prop:thickness-La(x)}, Lemma \ref{lem:thickness-interval-sequences} and Lemma \ref{lem:interval-comparison} there   exist $N\in\N$ such that for any $i,j>N$ we have
  \begin{equation}
    \label{eq:14-3}
    \tau(E_i(x))>\eta,\quad \tau(E_j(y))>\eta;
  \end{equation}
  \begin{equation}\label{eq:13-6}
  \theta_i(x)=\min\set{\frac{|I_i^x|}{|G_i^x|}, \frac{|I_{i+1}^x|}{|G_i^x|}}>\eta, \quad \theta_j(y)=\min\set{\frac{|I_j^y|}{|G_j^y|}, \frac{|I_{j+1}^y|}{|G_j^y|}}>\eta;
  \end{equation}
  and
   \begin{equation}
    \label{eq:15-4}
    \frac{|I_{i+1}^x|}{|I_i^x|}<\eta(x),\quad \frac{|I_{j+1}^y|}{|I_j^y|}<\eta(y).
  \end{equation}
  Then (\ref{eq:13-6}) suggests that for each  $i>N$ the length of  the gap  $G_i^x$ is strictly smaller than the length of each of its neighboring intervals $I_i^x$ and $I_{i+1}^x$. Also, for each $j>N$ the length of the gap $G_j^y$ is strictly smaller than {the length of} each of  its neighboring intervals $I_j^y$ and $I_{j+1}^y$.

  \emph{Claim 1:    there must exist a pair $(i,j)\in\N \times\N $ with $i,j>N$ such that $I_i^x\cap I_j^y\ne\emptyset$.}

 \noindent{\bf Proof of Claim 1.} Take $\de>0$ sufficiently small such that there exists a smallest integer $N_1>N$ such that $I_{N_1}^x\cap(1/m-\de, 1/m)\ne\emptyset$, and there exists a  smallest integer  $N_2>N$ such that $I_{N_2}^y\cap(1/m-\de, 1/m)\ne\emptyset$. Then $I_i^x\subset(1/m-\de, 1/m)$ for any $i>N_1$, and $I_j^y\subset(1/m-\de, 1/m)$ for any $j>N_2$. So to prove Claim 1 it suffices to prove that there  exist $i> N_1$ and $j>N_2$ such that $I_i^x\cap I_j^y\ne\emptyset$. Suppose this is not true. Then each basic interval $I_i^x$ with $i> N_1$ must belong to a gap of the interval sequence $\set{I_j^y: j\ge N_2}$. In other words, for each $i>N_1$ there must exists $j\ge N_2$ such that $I_i^x\subset G_j^y$. Observe that $p_i^x\nearrow1/m$ as $i\to\f$, and $p_j^y\nearrow 1/m$ as $j\to\f$. So there must exist $i> N_1$ such that $I_i^x$ and $I_{i+1}^x$ belong to two different gaps of the interval sequence $\set{I_j^y: j\ge N_2}$, say $I_i^x\subset G_j^y$ and $I_{i+1}^x\subset G_{j'}^y$ for some $j'>j\ge N_2$ (see Figure \ref{fig:3}).
    \begin{figure}[h!]
\begin{center}
\begin{tikzpicture}[
    scale=10,
    axis/.style={very thick, ->},
    important line/.style={thick},
    dashed line/.style={dashed, thin},
    pile/.style={thick, ->, >=stealth', shorten <=2pt, shorten
    >=2pt},
    every node/.style={color=black}
    ]

\draw[important line,blue] (-0.2, 0.7)--(-0.02, 0.7); \node[] at(-0.11,0.65){$I_{j}^y$};
    \draw[red,important line] (0, 0.7)--(0.4, 0.7);
     \node[] at(0.2, 0.75){$I_i^x$};
    \draw[red,important line] (1, 0.7)--(1.2, 0.7);
          \node[] at(1.1, 0.75){$I_{i+1}^x$};
\draw[{blue, important line}] (0.46, 0.7)--(0.6, 0.7); \node[] at(0.53, 0.65){$I_{j+1}^y$};
          \draw[dashed line, blue] (0.7,0.7)--(0.8,0.7);
      \draw[important line,blue] (0.9, 0.7)--(0.98, 0.7); \node[] at(0.94, 0.65){$I_{j'}^y$};
 \draw[important line,blue] (1.3, 0.7)--(1.36, 0.7); \node[] at(1.33, 0.65){$I_{j'+1}^y$};

\end{tikzpicture}
\end{center}
\caption{The geometrical structure of the closed intervals $I_i^x=[p_i^x, q_i^x]=conv(E_i(x))$ and  $I_{i+1}^x=[p_{i+1}^x, q_{i+1}^x]=conv(E_{i+1}(x))$. They belong  to the gaps $G_j^y$ and $G_{j'
}^y$ respectively.}\label{fig:3}
\end{figure}
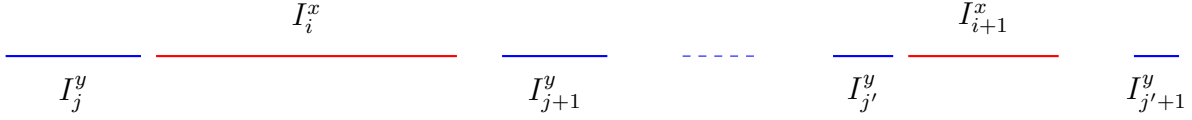
  Since $\theta_i(x)>1$ by (\ref{eq:13-6}), we have
 \[
 |I_{j'}^y|\le |G_i^x|<\min\set{|I_i^x|, |I_{i+1}^x|}.
 \]
This implies that
\[
\theta_{j'}(y)\le\frac{|I_{j'}^y|}{|G_{j'}^y|}\le \frac{|I_{j'}^y|}{|I_{i+1}^x|}<1,
\]
which leads to a contradiction with (\ref{eq:13-6}) that $\theta_j(y)>\eta>1$ for all $j>N$.  Hence, there must exist $i,j>N$ such that $I_i^x\cap I_j^y\ne\emptyset$, proving Claim 1.

\emph{Claim 2: there exists a pair $(i_1, j_1)\in\N\times\N$ with $i_1,j_1>N$ such that the two Cantor sets $E_{i_1}(x)$ and $E_{j_1}(y)$ are interleaved.}

 {\bf Proof of Claim 2.} By Claim 1 we can choose   $i, j>N$ such that $I_i^x\cap I_j^y\ne\emptyset$. If $I_i^x\nsubseteq I_j^y$ and $I_i^x\nsupseteq I_j^y$, then the two Cantor sets $E_i(x)$ and $E_j(y)$ are interleaved, and then we are done by setting $i_1=i$ and $j_1=j$. Otherwise, without loss of generality we may assume $I_i^x\subset I_j^y=conv(E_j(y))$. If $I_i^x$ is not contained in a gap of $E_j(y)$, then $E_i(x)$ and $E_j(y)$ are interleaved Cantor sets, and {again} we are done by setting $i_1=i$ and $j_1=j$. In the following we assume that $I_i^x$ is contained in a gap of $E_j(y)$, say
$I_i^x\subset G_j^y(\w_1)$
 for some word $\w_1$ of smallest length. Here $G_j^y(\w_1)$ is the gap between the two basic intervals $I_j^y(\w_1^+)$ and $I_j^y(\w_1)$ which are defined as in (\ref{eq:thickness-j}), see Figure \ref{fig:4}.
  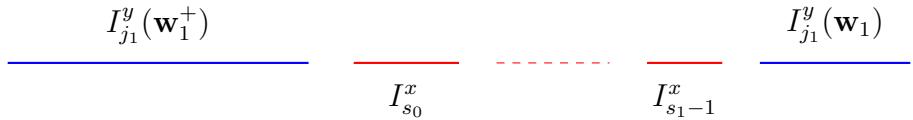
\begin{figure}[h!]
\begin{center}
\begin{tikzpicture}[
    scale=10,
    axis/.style={very thick, ->},
    important line/.style={thick},
    dashed line/.style={dashed, thin},
    pile/.style={thick, ->, >=stealth', shorten <=2pt, shorten
    >=2pt},
    every node/.style={color=black}
    ]

    \draw[blue,important line] (0, 0.7)--(0.4, 0.7);
     \node[] at(0.2, 0.75){$I_{j_1}^y(\w_1^+)$};
    \draw[blue,important line] (1, 0.7)--(1.2, 0.7);
          \node[] at(1.1, 0.75){$I_{j_1}^y(\w_1)$};
\draw[important line,red] (0.46, 0.7)--(0.6, 0.7); \node[] at(0.53, 0.65){$I_{s_0}^x$};
          \draw[dashed line, red] (0.65,0.7)--(0.8,0.7);

          \draw[important line,red] (0.85, 0.7)--(0.95, 0.7); \node[] at(0.9, 0.65){$I_{s_1-1}^x$};

\end{tikzpicture}
\end{center}
\caption{The geometrical structure of the closed intervals $I_{j_1}^y(\w_1^+), I_{s_0}^x, I_{s_1-1}^x$ and  $I_{j_1}^y(\w_1)$. The gap $G_{j_1}^y(\w_1)$ is the open interval between the two basic intervals $I_{j_1}^y(\w_1^+)$ and $I_{j_1}^y(\w_1)$.}\label{fig:4}
\end{figure}
We set $s_0=i$ and $j_1=j$. Note by (\ref{eq:14-3}) that the thickness of $E_{j_1}(y)$ is strictly larger than one. Then the  length of the interval $I_{j_1}^{y}(\w_1)$  is strictly larger than {the length of} the gap $G_{j_1}^y(\w_1)$. Furthermore, note by (\ref{eq:13-6}) that  the thickness of the interval sequence $\set{I_i^x: i>N}$ is larger than one, and $p_i^x\nearrow 1/m$ as $i\to\f$, we claim that there must exist a smallest integer $s_1>s_0$ such that
\begin{equation}\label{eq:14-1}
I_{s_1}^x\cap I^y_{j_1}(\w_1)\ne\emptyset.
\end{equation}
Suppose on the contrary that there is no such $s_1>s_0$ satisfying $I_{s_1}^x\cap I^y_{j_1}(\w_1)\ne\emptyset$. Then there is an $s>s_0$ such that $I_{j_1}^y(\w_1)$ is contained in the gap $G_s^x$ between $I_{s}^x$ and $I_{s+1}^x$. This implies that
\[
\theta_s(x)\le \frac{|I_s^x|}{|G_s^x|}\le \frac{|I_s^x|}{|I_{j_1}^y(\w_1)|}<1,
\]
where the last inequality follows by that $I_s^x\subset G_{j_1}^y(\w_1)$ and $|G_{j_1}^y(\w_1)|<|I_{j_1}^y(\w_1)|$. However, this leads to a contradiction with (\ref{eq:13-6}) that $\theta_i(x)>1$
for all $i>N$. Hence, we prove the existence of $s_1$  such that (\ref{eq:14-1}) holds.

Note by (\ref{eq:14-3}) and (\ref{eq:15-4})   that
\begin{equation}
\label{eq:15-5}|I_{s_1}^x|\le \eta(x)|I_{s_1-1}^x|\le \eta|G_{j_1}^y(\w_1)|<|I_{j_1}^y(\w_1)|.
\end{equation}
Since $I_{s_1}^x\cap I_{j_1}^y(\w_1)\ne\emptyset$,
  by (\ref{eq:15-5}) it follows that
 if $I_{s_1}^x$ is not contained in a gap of $E_{j_1}(y)\cap I_{j_1}^y(\w_1)$, then $E_{s_1}(x)$ and $E_{j_1}(y)$ are interleaved, and then we are done by setting $i_1=s_1$. Otherwise, suppose $I_{s_1}^x$ is contained in a gap of $E_{j_1}(y)\cap I_{j_1}^y(\w_1)$. Then there exists a word  $\w_2$ of smallest length such that $I_{s_1}^x\subset G_{j_1}^y(\w_1\w_2)$. Note that the thickness of $E_{j_1}(y)$ and the thickness of $\set{I_i^x: i\ge N}$ are both larger than one. Then by the same argument as in (\ref{eq:14-1}) there exists a smallest integer  $s_2>s_1$ such that
\[
I_{s_2}^x\cap I^y_{j_1}(\w_1\w_2)\ne\emptyset.
\]
Again, by the same reason as in the proof of (\ref{eq:15-5}) we can prove that $|I_{s_2}^x|<|I_{j_1}^y(\w_1\w_2)|$. So, if $I_{s_2}^x$ is not contained in a gap of $E_{j_1}(y)\cap I_{j_1}^y(\w_1\w_2)$, then $E_{s_2}(x)$ and $E_{j_1}(y)$ are interleaved, and then we set $i_1=s_2$. Otherwise, there exists a word $\w_3$ of smallest length such that $I_{s_2}^x\subset G_{j_1}^y(\w_1\w_2\w_3)$.

Repeating the above argument, and we claim that our procedure must stop at some finite time $n$, i.e.,  there   exist  an integer   $s_n$    and  words $\w_1, \w_2,\ldots, \w_n$ such that
\[I_{s_n}^x \cap I_{j_1}^y(\w_1 \w_2 \cdots\w_n)\ne\emptyset, \quad |I_{s_n}^x|<|I_{j_1}^y(\w_1\w_2\ldots \w_n)|,\] and  $I_{s_n}^x$ is not contained in a gap of  $E_j(y)\cap I_{j_1}^y(\w_1 \w_2 \cdots\w_n)$. If {this} claim is true, then we are done by setting $i_1=s_n$. Suppose {it} is not true. Then  $I_i^x\subset I_{j_1}^y$ for all large integers $i>s_1$, which implies that  $q_i^x\le q_{j_1}^y<1/m$ for all $i\ge 1$. This leads to a contradiction with $\lim_{i\to\f}q_i^x= 1/m$.
Therefore, we have find a pair $(i_1,j_1)\in\N\times\N$ with $i_1, j_1>N$ such that the two Cantor sets $E_{i_1}(x)$ and $E_{j_1}(y)$ are interleaved. This proves Claim 2.

\medskip

Now we take $\de>0$ small enough such that $1/m-\de> {\max}\set{q_{i_1}^x, q_{j_1}^y}$, and repeating the   argument as in Claim 1 and Claim 2. Then we can find another pair $(i_2, j_2)\in\N\times\N$ with $i_2>i_1, j_2>j_1$ such that $E_{i_2}(x)$ and $E_{j_2}(y)$ are interleaved. We can proceed this argument  indefinitely, and then we find a  sequence of pairs $(i_k, j_k)\in\N\times\N$ such that for any $k\ge 1$ we have  $i_k<i_{k+1}, j_k<j_{k+1}$, and the two Cantor sets $E_{i_k}(x)$ and $E_{j_k}(y)$ are interleaved. This completes the proof.
\end{proof}

\begin{proof}
  [Proof of Theorem \ref{th:intersection-La(x)-La(y)}]
  By Lemmas \ref{lem:thickness-intersection} and \ref{lem:interleaved-pairs} it follows that $\La(x)\cap \La(y)$ contains a Cantor set of thickness arbitrarily large. More precisely, for the infinite sequence of pairs $(i_k, j_k) \in\N\times\N$ defined as in Lemma \ref{lem:interleaved-pairs} it follows that  for any $k\ge 1$,
  \begin{equation}\label{eq:14-2}
  \La(x)\cap\La(y)\supset E_{i_k}(x)\cap E_{j_k}(y)
  \end{equation}
and the intersection  $E_{i_k}(x)\cap E_{j_k}(y)$ contains a Cantor  set with its thickness at least of order $\sqrt{\min\set{\tau_{\mathcal N}(E_{i_k}(x)), \tau_{\mathcal N}(E_{j_k}(y))}}$. Note by Lemma \ref{lem:thickness-two-defs} and Proposition \ref{prop:thickness-La(x)} that \[\lim_{k\to\f}\tau_{\mathcal N}(E_{i_k}(x))=\lim_{k\to\f}\tau_{\mathcal N}(E_{j_k}(y))=+\f.\] Then by (\ref{eq:dim-thickness}) and (\ref{eq:14-2}) we conclude that
\[
\dim_H(\La(x)\cap\La(y))\ge \dim_H(E_{i_k}(x)\cap E_{j_k}(y))\to \; 1 \quad\textrm{as }k\to\f.
\]
This completes the proof.
\end{proof}

{\section*{Acknowledgements}
 The first author was supported by NSFC No.~11701302,  Zhejiang Provincial Natural Science Foundation of China with No.LY20A010009, and the K.C. Wong Magna Fund in Ningbo University. The second author was supported by NSFC No.~11971079 and the Fundamental
and Frontier Research Project of Chongqing No.~cstc2019jcyj-msxmX0338 and No.~cx2019067. The third author was  supported by NSFC No.~12071148 and Science and Technology Commission of Shanghai Municipality (STCSM)  No.~18dz2271000.}



\end{document}